
\def\real{{\tt I\kern-.2em{R}}}  
\def\nat{{\tt I\kern-.2em{N}}}   
\def\eps{\epsilon}
\def\realp#1{{\tt I\kern-.2em{R}}^#1}
\def\natp#1{{\tt I\kern-.2em{N}}^#1}
\def\hyper#1{\,^*\kern-.2em{#1}}
\def\monad#1{\mu (#1)}

\def\st#1{{\tt st}(#1)}
\def\hyperreal{{^*{\real}}}
\def\hyperrealp#1{{\tt ^*{I\kern-.2em{R}}}^#1} 
\def\hypernat{{^*{\nat }}}
\def\hypernatp#1{{{^*{{\tt I\kern-.2em{N}}}}}^#1} 
\def\eskip{\hskip.25em\relax}

\def\Hyper#1{\hyper {\eskip #1}}
\def\leaderfill{\leaders\hbox to 1em{\hss.\hss}\hfill}
\def\srealp#1{{\rm I\kern-.2em{R}}^#1}

\def\pars{\par\smallskip}
\def\parm{\par\medskip}

\def\b#1{{\bf #1}}
\def\ref#1{$^{#1}$}

\def\sig{{^\sigma}}
\def\m@th{\mathsurround=0pt}
\def\rightarrowfill{$\m@th \mathord- \mkern-6mu \cleaders\hbox{$\mkern-2mu 
\mathord- \mkern-2mu$}\hfil \mkern-6mu \mathord\rightarrow$}
\def\leftarrowfill{$\mathord\leftarrow
\mkern -6mu \m@th \mathord- \mkern-6mu \cleaders\hbox{$\mkern-2mu 
\mathord- \mkern-2mu$}\hfil $}
\def\noarrowfill{$\m@th \mathord- \mkern-6mu \cleaders\hbox{$\mkern-2mu 
\mathord- \mkern-2mu$}\hfil$}
\def\orgate{$\bigcirc \kern-.80em \lor$}
\def\andgate{$\bigcirc \kern-.80em \land$}
\def\inverter{$\bigcirc \kern-.80em \neg$}
 \magnification=1200
\tolerance 10000
\hoffset=0.25in
\hsize 6.00 true in
\vsize 8.75 true in

\def\id{\par\hangindent2\parindent\textindent}
\def\textindent#1{\indent\llap{#1}}
\font\eightrm=cmr9
\centerline{\bf Nonstandard Analysis and Generalized Functions}\bigskip
\centerline{Robert A. Herrmann}\par\medskip
\centerline{Mathematics Department}
\centerline{U. S. Naval Academy}
\centerline{572C Holloway Rd.}
\centerline{Annapolis,  MD 21402-5002}
\centerline{20 NOV 1994}
\centerline{\ }\bigskip
{\leftskip=0.5in \rightskip=0.5in \noindent {\eightrm {\it Abstract:} This application of nonstandard analysis utilizes the notion of the highly-saturated enlargement. These nonstandard methods are applied to the theory of generalized functions (distributions) and demonstrates how such analysis clarifies many aspects of this theory. \par}}\par
\vskip 14pt
\leftline{\bf 1. Additional Modeling Concepts.}\par    
\medskip 
In what follows$,$ the basic notation and definitions are as they  appear in 
{[3]}. However$,$ 
except as mentioned in the appendix of [3]$,$ a significant aspect of 
nonstandard analysis has not been  developed fully. The 
superstructure  constructed in [3] is a 
model for $\Gamma$ where
$\Gamma$ is the set of all sentences that hold in the structure. 
Since it is constructed from $\real$$,$ it is called a model (for real     
analysis)$,$ for apparently every true statement from analysis   
holds true in the superstructure.  
The  
 elementary nonstandard structure  
$\hyper {\cal M} = (\Hyper {\cal H}, \in , = )$ is associated with the 
standard model ${\cal M} = ({\cal H}, \in , = )$ in a slightly special sense. 
Due to the applications in [3]$,$ it was not necessary to discuss $\hyper {\cal 
M}$ relative to its special properties. This is no longer the case.\par
The structure $\hyper {\cal M}$ is constructed from a (bounded) ultrapower  
based upon the structure $\cal M$. {[1$,$ p. 15--19]}$,$ {[2]$,$} {[5$,$ p. 83--88]} 
It is assumed that there are 
constants 
that denote every member of $\cal H$ where we do not differentiate between a 
constant and the object it names. Suppose that $J$ is the index set and that 
$\cal U$ is an appropriate ultrafilter on $J.$ Let infinite $A\in \cal H$ and 
$A$ contains no individuals in $\real.$ 
Let 
$f \in {\cal H}^J$ and $f(j) = A,$ for each $j \in J.$ This is the constant 
map$,$ constant in two ways both as a mathematical entity and relative to the 
value being denoted by a constant in our language. Some authors define an 
injection $k$$,$ which is denoted by  $e$ in [5] and $i$ in [1]$,$ such that 
$k(A) = [A],$ where $[A]$ is the $\cal U$-equivalence class in the ultrapower
that contains the constant  map $f.$ To obtain an isomorphic copy of $\cal 
M,$  we could follow the usual Mostowski collapsing process as 
described in [5$,$ p. 84--85] that gives the $\hyper {}$ mapping and let  
this map be restricted to the domain of all the constant sequence $\cal U$-equivalence 
classes in the ultrapower. This leads to an isomorphic copy 
of $\cal M.$ [5$,$ p 85.] But$,$ care must be 
taken relative to the interpretation of the * map. It must always be 
remembered the $\hyper {\cal M}$ is a model of the bounded expressions that 
hold in $\cal M$ although sometimes the bounding set is not expressly stated 
in an expression it must be understood that the quantifiers are restricted to 
specific elements in some $X_n$ for $\cal M$ and for the corresponding 
$\Hyper X_n$ for
 $\Hyper {\cal H}$. The members of $\Hyper X_n$ are ``internal'' entities.  \par
Consider the first-order statement $a \in b$ about objects in $\cal H.$
Then the statement $\hyper a \in \Hyper b$ holds in $\hyper {\cal M}$ and 
is relative to the embedding of $\cal M$ into $\hyper {\cal M}.$
Now consider a set $A \in \cal H$ and the bounded statement 
$S = \forall x ((x \in A)\land(x \in b)),$ which is in the required form, and let $S' = \{x\mid (x \in A)\land(x \in b)\}.$ 
Then $S$ holds in $\cal M$  if and only if    
$\Hyper S = \forall x ((x \in \hyper A)\land (x \in \Hyper b))$ 
holds in $\hyper {\cal M}.$
This indicates that the 
quantifiers are restricted to members of $\hyper A.$ Since $A$ is a set,
$A \in X_n$ for some $n \geq 1.$ Then from Proposition 1 (iv) [3]$,$
it follows that each specific ``$x$'' that satisfies $\Hyper {S'}$ is a member of 
$\Hyper X_0$ or $\Hyper X_n,$ where $n\geq 1.$ Thus each such ``$x$'' is a 
member of $\Hyper {\cal H}.$ \par
On the other hand$,$ if $S$ holds in $\cal M,$ then $S$ also holds in the 
isomorphic copy of $\cal M,$  where the quantifier is restricted to 
constant sequence ${\cal U}$-equivalence classes. Hence the only members in the isomorphic copy of 
the set $S'$ are the $\hyper a$ 
such that $(\hyper a \in \hyper A)\land(\hyper a \in \Hyper b).$ But is this isomorphic copy of $S'$ a 
member of  
$\Hyper {\cal H}$? One of the ideas behind the concept of a nonstandard 
structure  
is that 
if $S'$ is infinite$,$ then its isomorphic copy is not a member of $\Hyper {\cal 
H}$ although it is a subset of $\hyper A$ and $\hyper A \in \Hyper {\cal H}$. 
In this case$,$ the set $\Hyper {S'}$ contains entities that are not produced by 
the * map. In order to discuss sets such as the isomorphic copy of $S',$
a second superstructure $({\cal Y},\in,=)$ is generated with ground set $Y_0 = \hyper X_0.$ 
It is within 
this superstructure that we can specifically construct the sets that restrict 
the quantifiers when bounded statements are interpreted for the isomorphic 
copy of $\cal M.$ \par
\vskip 14pt
\hrule  
\smallskip
{\bf Definition 1.1 {(Sigma Operator$,$ Standard Copy)}} If a set 
$A \in \cal H,$ then $\sig A = \{\hyper a\mid a \in A \},$ where the $\hyper  a$ in this form denotes the constant sequence $\cal U$-equivalence class [8, p. 44-45].  
 \par
\smallskip
\hrule
\vskip 14pt
\noindent Let $\sig {\cal H} = \{\sig A\mid A \in 
{\cal H}\}.$  There is a certain 
confusion of symbols that one tries to avoid. When the isomorphic embedding is 
being considered$,$ it is understood that $\hyper a$ means something 
different when viewed as a {\it set} of objects. The symbol $\hyper a$ is used 
only as a name for the equivalence class in the ultrapower that contains the
constant sequence for the isomorphic embedding. But, relative to the structure $\hyper {\cal M}$$,$ and as a nonempty set, $\hyper a$ 
contains the objects in $\sig a.$ Thus relative to the sets in $({\cal Y},\in 
=),$  the isomorphic copy is determined by the $\sig$ and $\hyper {}$ operators. It is this 
fact that has lead to a minimizing of the use of the $\sig$ notation. 
Thus it is customary to do our real analysis
in $\cal M$ rather then in $\sig {\cal M} = (\sig {\cal H},\in,=)$ knowing that$,$ when {\it 
comparisons} are 
to be made$,$ we can apply the isomorphism to obtain the actual objects in
$\sig {\cal H}$ that are used for  modeling purposes. I mention that significant fact is that if a set $A \in \cal H$ is finite$,$ then
$\hyper A = \sig A.$ One might say that the map $\sig$ preserves the 
hierarchy of finite sets. (Some authors find no need to consider $\sig {\cal M}$ [5].) \par

From the viewpoint  of abstract model theory$,$ what are the real numbers?
The real numbers is considered to be ANY structure isomorphic to the 
``standard'' structure. The standard structure is considered to be a nonempty 
set $\real$ with the appropriate operators and unique elements defined and 
such that the operators$,$ unique elements and subsets of $\real$ satisfy a set 
of axioms. Under the above isomorphism *$,$ Theorem 3.1.1 in [3] implies that the 
isomorphic copy of the real numbers can be considered as THE real numbers
and real analysis takes place in $\sig \cal  H.$
Thus$,$ as has become customary$,$ we let $\sig \real = \real.$ Moreover,
to more fully express this identification of the real numbers$,$ 
consider how this isomorphism deals with the 
operator $\subset.$  
From the definition of the operator $\subset,$ it follows that if
$A \in \real,$ then $\sig A = \emptyset.$ Further$,$ given two 
$A, B \in {\cal H}.$ Then $ A\subset B$ if and only if $\sig A\subset\ 
{\sig B}.$ Hence$,$ if $A \subset \real,$ then $A =\ {\sig A}\subset \hyper 
A.$ Basic operators such as $+,\  \cdot,\ <,$ under the isomorphism$,$ 
become the 
operators $\sig +,\  \sig \cdot,\ \sig <$ which then become THE 
operators $+,\  \cdot,\ <$ for the field $\real.$ Although the $\sig$ 
notation could continue to be used on sets at any point when one is in doubt$,$ 
it has become customary to remove this notation in some cases like  $\sig \real  =\real.$  
Often relations such as $\hyper +,\ \hyper{\cdot},\ 
\hyper <$ for the ordered field $\hyperreal$  are considered as the actual 
extensions of the relations   $+,\ {\cdot},\ <.$ \par
Relative to notation$,$ what this means is that many of the constants
in $C({\cal H})$ that denote the members of $\cal H$ will now be used
to denote many members of $\sig \cal H.$ New symbols $\sig A$ are used 
to denote other members of $\sig {\cal H}.$  Then we have objects in 
$\hyper {\cal H}$ that have names in the symbol set $C(\hyper {\cal H})$ and
represent the internal extended standard objects. These are denoted by 
use of the * notation. All other members of $\hyper {\cal H},$ other than the 
internal members$,$ are denoted by distinctly different symbols
in $C(\hyper {\cal H}).$  Since there are only so many symbols that can be 
used$,$ we must state that this non-``started''  symbol represents an ``internal'' object.
Each object in $({\cal Y},\in,=)$ has a symbol name. All such objects 
that are not denoted by any previously defined notation must be explicitly
defined as ``external'' objects. In this regard$,$ if infinite $A \in {\cal H},$
then for the actual structure $\hyper {\cal M}$ soon to be 
constructed$,$ the object (denoted by) $\sig A$ is external.  \par
Relative to mathematical modeling$,$ if we are modeling physical entities with 
names taken from a specific discipline dictionary$,$ then it is immaterial which 
real analysis structure is used for the modeling correspondence. Thus$,$ the
isomorphic $\sig {\cal M}$ is chosen as the 
appropriate structure. Hence$,$ if the Euclidean n-space 
function $\b v(t)$ corresponds to a natural physical
world (i.e. N-world) vector$,$ then the function $\Hyper {\b v(t)}$ corresponds to a nonstandard 
physical world (i.e. NSP-world) vector. Since $\b v \subset \Hyper {\b 
v},$ and $\b v,\Hyper {\b v} \in ({\cal Y}, \in, =),$ the N-world physical vector can be assumed to be a restriction of the 
NSP-world vector to the N-world. Although more can be said about the effects 
of such a restriction relative to direct and indirect observation$,$ it is not 
necessary$,$ at this point$,$ to delve more deeply into such concepts.\par 
There has also arisen a certain terminology. Suppose that $f$ is a continuous 
map from $\real$ into $\real$ (i.e. $f \colon \real \to \real$). Then one can 
write that $\hyper f$ is a *-continuous map from $\hyperreal$ to $\hyperreal.$
The term  ``*-continuous'' is often  replaced in scientific discourse by the 
term   ``hypercontinuous.''  On the other hand$,$ some authors leave the term in 
this form and when communicating  orally say ``hypercontinuous'' or
``star-continuous.'' \par
\vskip 14pt\vfil\eject
\hrule  
\smallskip
{\bf Definition 1.2 {(Concurrent Relation)}} A (bounded) binary relation $\Phi$
in  $\cal H$ and$,$ hence$,$ in $\sig {\cal H},$ is {\it concurrent} if the 
following holds. For each finite $\not= \emptyset$ set 
$A = \{(a_1,b_1),\ldots,(a_n,b_n)\} \subset \Phi$ there exists in the range$,$ 
$R(\Phi),$ of $\Phi$ 
some $b$ such that $\{(a_1,b),\ldots,(a_n,b)\} \subset \Phi.$ \par
\smallskip
\hrule
\vskip 14pt
Now as to the actual construction of the nonstandard structure $\hyper {\cal 
M},$  a special ultrafilter is selected which has the following  {\it 
enlarging} property. \par 
\vskip 14pt
\hrule  
\smallskip
{\bf Definition 1.3 {(Enlargement)}} The structure $\hyper {\cal M}$ is an 
{\it enlargement} if for every concurrent relation $\Phi$  with the domain 
$D(\Phi)$ there exists an 
internal $b\in \Hyper {(R(\Phi))} = R(\Hyper {\Phi)}$ such that for each 
$\hyper a \in \sig D(\Phi),\ (\hyper a, b) \in \Hyper {\Phi}.$   
\par
\smallskip
\hrule
\vskip 14pt
\noindent In all that follows$,$ it is assumed that $\hyper {\cal M}$ is$,$ at the least$,$ an 
enlargement.  Also note that all members of $\Hyper {(R(\Phi))}$ are 
internal.\par
\vskip 14pt
{\bf Theorem 1.1.} {\sl Let infinite $A \in \cal H.$ Then $\sig A \subset 
\hyper A$ and $\sig A \not= \hyper A.$}\pars
Proof. Suppose that infinite $A \in \cal H.$ Consider the relation
$\Phi = \{(x,y)\mid (x \in A)\land(y\in B)\land(x \not= y)\}.$ Suppose that 
$\{(a_1,b_1),\ldots,(a_n,b_n)\} \subset \Phi.$ However$,$ since $A$ is 
infinite$,$ there exists some $b \in A$ such that $b \not=a_i,\ i = 1,\ldots, 
n.$ Hence$,$ $\{(a_1,b),\ldots,(a_n,b)\} \subset \Phi.$ Thus $\Phi$ is a 
concurrent relation.  Consequently$,$ there is an internal $b \in R(\Hyper {\Phi})$ such 
that $(\hyper a, b) \in \Hyper {\Phi}$ and $b \not= \hyper a.$ This 
completes the proof.\par
\vskip 14pt
Often just one identified concurrent relation can determine a major portion of 
an entire nonstandard theory. A few more examples indicate this fact. 
First$,$ consider the extension of the absolute value function to $\hyperreal.$
By definition$,$ for any $r \in \real$ if $r \geq 0,$ then $\vert r \vert  = r$
and for any $r \in \real$ if $r < 0,$ then $\vert r \vert  = -r.$ Stated 
formally$,$ we have that $S = \forall x(((x\in \real)\land (r \geq 0))\to
(\vert r \vert = r))\land ((x\in \real)\land (r < 0))\to
(\vert r \vert = -r))).$ This statement holds in $\cal M.$  Thus its 
*-transform holds in $\hyper {\cal M}.$ The *-transform is 
$\Hyper S = \forall x(((x\in \hyperreal)\land (r \geq 0))\to
(\Hyper {\vert r \vert)} = r)\land ((x\in \hyperreal)\land (r < 0))\to
(\Hyper {\vert r \vert)} = -r)),$  where $\vert \cdot \vert$ is viewed as a unary
operator. Thus the operator $\Hyper {\vert \cdot \vert}$ is but the absolute 
value operator as it is defined for the totally ordered field $\hyperreal.$ 
Hence we can drop the * notation from $\Hyper {\vert \cdot \vert}.$
\par
\vskip 14pt
{\bf Theorem 1.2.} {\sl There exists in $\hyperreal,$ a nonzero infinitesimal.} 
\par
\smallskip Proof. Consider the relation $\Phi = \{(n,m) \mid (0<(1/m)<(1/n))\land(n \in 
\nat)\land(m \in \nat)\}.$ Suppose that $\{(n_1,m_1),\ldots,(n_j,m_j)\} 
\subset \Phi.$ Let $M= \max \{m_1,\cdots,m_j\}.$ Then 
$\{(n_1,M+1),\ldots,(n_j,M+1)\} \subset \Phi$ implies that $\Phi$  is 
concurrent. Hence$,$ there exists some $\Lambda \in \hypernat$ such that 
$ 0<1/\Lambda < 1/n$ for all $n \in \sig\nat = \nat.$ Now consider any 
positive 
$r \in \real.$ Then there exists some $n \in \nat$ such that $0<1/n<r.$ Hence,
$\vert 1/\Lambda \vert < r.$ But since $r$ is an arbitrary positive real 
number$,$ this last results hold for all $r \in \real^+.$ This completes the 
proof.\par
\vskip 14pt
An examination of chapter 2 in [3] shows that the properties of the set of 
all infinitesimals $\monad 0$ are determined from Theorem 1.2. One of the most 
significant portions of the nonstandard theory of analysis is relative to the 
set of all  ``hyperfinite'' sets. These are all of the internal sets $A\in 
\Hyper (F({\cal H})= \cup \{\Hyper {(F(X_n))}\mid X_n\in {\cal H}\}.$ 
Also note that if $A \in X_n,\ (n >0),$  then $F(A) \in X_{n+1}.$  Hence,
$\Hyper {(F(A))} \in \hyper {X_{n+1}}.$ Viewed as a mapping$,$ we have that 
$\Hyper F$ is defined on all internal sets in $\Hyper {\cal H}$  and if 
$B \in \cal H$$,$ then $\Hyper {(F(B))} = \hyper F(\hyper B).$ As shown in [3] Theorem 
4.3.2$,$ the above definition for the hyperfinite sets is equivalent to the 
the internal bijection definition. 
The hyperfinite sets satisfy in $\hyper {\cal M}$ all of the 
first-order properties associated with finite sets. But from the exterior 
viewpoint of nonstandard analysis$,$ such sets are far from being finite.\par
\vskip 14pt
To get the full strength of nonstandard analysis as it relates to generalized 
functions$,$ we need to concept of the $\kappa$-saturated enlargement$,$ where 
$\kappa$ is a infinite cardinal number. \par
\vskip 14pt
\hrule  
\smallskip
{\bf Definition 1.4 {( $\kappa$-Saturated )}} 
The structure $\hyper {\cal M}$ is a {\it $\kappa$-saturated} 
if given any internal (bounded) 
binary relation $\Phi,$  with the internal domain $D(\Phi),$ that is 
concurrent on  
$A \subset D(\Phi),$  where cardinality of $A< \kappa,$ then there exists 
 an 
internal $b\in R(\Phi),$ the internal range of $\Phi,$ such that for each 
$a \in A,\ (a, b) \in {\Phi}.$ \par
\smallskip
\hrule
\vskip 14pt
Throughout the remaining portions of this paper$,$ we assume that 
$\hyper {\cal M}$ is a $\kappa$-saturated$,$ where $\kappa$ is any 
(regular) cardinal number
greater than the cardinality of $\cal H.$ This would also imply that $\Hyper {\cal M}$
is an enlargement. It can be shown by means of the 
ultralimit process$,$ such  (bounded) $\kappa$-saturated enlargements exist. \par
\vskip 14pt
{\bf Theorem 1.3.} {\sl Consider internal $B$ and any $A \subset B$ such that 
cardinality of $A< \kappa.$ 
Then there exists a 
hyperfinite set $\Omega$  such that $A \subset \Omega\subset B.$}\par         
\smallskip Proof. From the construction of $\cal H$ we know that there is some $X_n,\ n 
\geq 1$ and 
$B \in \hyper {X_n},$  $\Hyper F(B) \in \hyper X_n,$ if $q \in B,\ \{q\}
 \in X_n$ and $q \in \hyper X_0 \cup \hyper X_{n-1}.$ 
Consider the internal binary relation $\Phi=\{(x,y)\mid (x \in y \in 
\Hyper F(B))\land(x \in \hyper X_0 \cup \hyper X_{n-1}) \land (y \in X_n)\}.$ 
Suppose that 
$\{(a_1,b_1),\ldots,(a_j,b_j)\} \subset \Phi.$ By *-transfer of the standard 
theorem$,$ it follows that $b = b_1 \cup \cdots \cup 
b_j \in \Hyper F(B)$ and $a_i \in b,\ i = 1,\ldots,j.$ Hence$,$ $\Phi$ is a 
concurrent on its domain. But $A \subset D(\Phi)$ and has the appropriate 
cardinality. Hence$,$ there exists some $\Omega\in 
\Hyper (F(B))$ such that $ A \subset \Omega\subset B.$ This completes the proof. \par
{\bf Corollary 1.3.1.} {\sl Consider standard $A.$ 
Then there exists a 
hyperfinite set $\Omega$  such that $\sig A \subset \Omega \subset \hyper 
A.$}\par          
\smallskip Proof. Simply note that the cardinality of $A$ is less than $\kappa.$\par
 \vskip 14pt
\noindent The hyperfinite sets are 
the basic building blocks of the nonstandard theory of probability spaces. 
\par \bigskip
 \noindent {\bf 2. Generalized Functions.}\par \medskip
 The 
functions considered are real valued functions. It is not difficult to 
extend all of the results in this section to complex valued functions. Further 
all standard functions map $\real$ into $\real.$ Let $C^\infty$ be the set of 
all real valued functions defined on $\real$ which have derivatives of all 
orders at each $x \in \real.$ The set $\Hyper C^\infty$  contains some very 
interesting *-continuous and *-differentiable functions. Throughout this 
paper$,$ nonempty ${\cal D}\subset C^\infty$ is always the notation for 
what is called the {\it test space.} Each member of ${\cal D}$ must be a 
function with bounded support. This implies that if $g \in {\cal D},$  then 
there is some $c \in \real$ such that $g(x) = 0$ for all $\vert x \vert \geq 
c.$  \par
 Usually one is interested in the generation of linear functionals. 
The customary generating functions are maps from $\real$ into $\real.$  The 
basic method of generation is by integration. Usually$,$ the customary 
integration is Lebesgue integration although it appears the generalized 
Riemann integral can also be used. The reason that Lebesgue is useful is that 
this integral applies to many highly discontinuous standard functions$,$ has 
useful convergence properties and$,$ operationally$,$ is sufficient. For our 
purposes$,$ the Lebesgue integral is considered as an operator in the sense 
that it is 3-tuple with the first coordinate a function$,$ the second an 
interval (or for other applications a measurable subset of $\real$)$,$ and the 
third coordinate the value when it exists. \par
 Our customary standard 
generating functions$,$ $CS$$,$ have the property that they are Lebesgue 
measurable on $[a,b],$ for $a \leq b, \ a,b \in \real$ and the integral 
$\int_a^b (f(x))^2 \, dx \in \real$ (i.e. $f \in {\cal L}^2([c,d])$ a 
classical Banach Space). If $f$ is measurable and bounded on $[a,b]$ and  $f 
\in {\cal L}([a,b])$ then $f \in {\cal L}^2([a,b]).$ {[7$,$ p. 219]} 
It is known that if $f \in {\cal L}^2(E),$ and for the Lebesgue measure$,$ $m$$,$ 
$m(E) \in \real,$ then $f \in {\cal L}(E).$ {[7$,$ p. 220]} 
 From this$,$ it 
follows that $\int_{-\infty}^\infty f(x)g(x)\, dx \in \real$ for each
$g \in {\cal D}.$ 
Our functions are restricted to members of internal set 
$\cap\{\hyper {\cal L}^2([c,d])\mid
(c\leq  d)\land(c,d \in \hyperreal)\}$  
so that the 
*-transform  of the classical Schwarz inequality applies.    
For the purposes of this paper$,$ the set of internal generating functions 
is has a slightly different formation and can only be assumed to be an 
external subset of $\cap\{\hyper {\cal L}^2([c,d])\mid
(c\leq  d)\land(c,d \in \hyperreal)\}.$ 
Recall 
that $\cal O$ is the set of all limited numbers in $\hyperreal.$ Note that 
this set is also called the set of finite numbers. It is clear that if
$f$ is a customary standard function$,$ then for each $c \leq  d,\  c, d \in 
\hyperreal,\  \hyper {\int_c^d}
(\hyper f(x))^2\, dx \in \hyperreal.$ Thus  $\hyper{\int_c^d}
(\hyper f(x))^2\, dx \in \hyperreal$  when $c,d \in {\cal O}.$ Further,
$\int_{-\infty}^\infty f(x)g(x)\, dx = r$ implies that $\hyper {\int_{-
\infty}^\infty} \hyper f(x)\hyper g(x)\, dx = \hyper r=r.$ First$,$ are there
members of $CS$ such that  $\hyper {\int_a^b}
(\hyper f(x))^2\, dx \in {\cal O}$  when $c,d \in {\cal O}$?\par
\vskip 14pt
{\bf Theorem 2.1.} {\sl Let $f\colon \real \to \real.$ Suppose that $f$ is 
measurable (in the Lebesgue sense) and bounded on all intervals 
$[a,b], \ a\leq b, \ a,b \in \real.$
Then for each $c\leq d, \ (c,d \in 
\hyperreal),\ \hyper {\int_c^d}
\hyper f(x)\, dx \in \hyperreal$ and $\hyper {\int_c^d}
(\hyper f(x))^2\, dx \in \hyperreal$ and if $c,d \in {\cal O},$ then
$\hyper {\int_c^d} \hyper f(x)\, dx \in {\cal O}$ and 
$\hyper {\int_c^d} (\hyper f(x))^2\, dx \in {\cal O}.$}
\par
\smallskip Proof. It follows from the hypotheses$,$ that  
$f \in {\cal L}([a,b])$ and for each interval $[a,b]$ there exists some $M\in 
\real$ such that $\vert f(x) \vert < M,$ for each $x \in [a,b].$ 
By *-transfer$,$ $ \hyper {\int_c^d}
\hyper f(x)\, dx \in \hyperreal$ for $c\leq d,\ (c,d \in\hyperreal).$ Let $c\leq d$ and $c,d \in {\cal O}.$ Then there exist $a,b \in 
\real$ such 
that $c \approx a,\ d \approx b.$ There are four cases to consider but one 
will suffice as a prototype.  Suppose that $c\leq a,\ d\leq b.$ Let $m$ denote 
the Lebesgue measure on the measurable subsets of $\real.$ For each $x, y \in 
\hyperreal$ such that $x \leq y,$ $\hyper m([x,y]) = y-x$ by *-transfer. Hence$,$ 
$\hyper m([c,a])  = c-a \approx 0,\ \hyper m([d,b]) = b-d \approx 0.$ 
Now $\hyper {\int_c^b} \hyper f\, dx =  \hyper {\int_c^a} \hyper f\, dx 
+\hyper {\int_a^d} 
\hyper f\, dx +\hyper {\int_d^b} \hyper f\, dx$ by *-transfer. Consider
$\hyper {\int_c^a} \hyper f\, dx.$ There exists some $g\in \real$ such that
$g\leq c.$ Hence $[c,b] \subset \Hyper [g,b].$ There exists some $M \in \real$ 
such that $\vert f(x) \vert \leq M$ for each $x \in [c,b].$ Again by 
*-transfer $\vert \hyper f(x) \vert \leq M$ for each $ x \in \Hyper [c,b].$ 
By *-transfer$,$ 
$$-M(\hyper m([c,a])) = -M(a-c) \leq \hyper {\int_c^a} \hyper f\, dx \leq
M(\hyper m([c,a]) = M(a-c).$$
Consequently$,$ $\hyper {\int_c^a} \hyper f\, dx \approx 0.$  In like manner.
$\hyper {\int_d^b} \hyper f\, dx \approx 0.$ Therefore$,$ $\hyper {\int_c^d} 
\hyper f\, dx
\approx \hyper {\int_a^b} \hyper f\, dx =r \in \real.$ Hence$,$ 
$\hyper {\int_c^d}
\hyper f(x)\, dx \in {\cal O}.$ \par
For the second part$,$ simply consider the known standard result that 
if $f \in {\cal L}([a,b])$ and $g$ is bounded and measurable on $[a,b],$ then 
$fg \in {\cal L}([a,b])$ and apply the first part. 
This completes the proof.\par
{\bf Corollary 2.1.1.} {\sl Let continuous $f\colon \real \to \real.$ 
Then for each $c\leq d, \ (c,d \in 
\hyperreal),\ \hyper {\int_c^d}
\hyper f(x)\, dx \in \hyperreal$ and $\hyper {\int_c^d}
(\hyper f(x))^2\, dx \in \hyperreal$ and if $c,d \in {\cal O},$ then
$\hyper {\int_c^d} \hyper f(x)\, dx \in {\cal O}$ and 
$\hyper {\int_c^d} (\hyper f(x))^2\, dx \in {\cal O}.$}\par
\smallskip Proof. Clearly$,$ $f\in {\cal L}([a,b])$ and $f\in {\cal L}^2([a,b]).$ The same 
proof as theorem 2.1 yields the conclusions.\par
\vskip 14pt
\hrule
\smallskip 
{\bf Definition 2.1 {(Generalized Functions)}} Let $T$ be the set of internal 
functions such that for each $f\in T,$ (i) $f\colon \hyperreal \to \hyperreal$ and
(ii) $\hyper {\int_c^d} (f(x))^2\, dx \in \hyperreal$ for each pair $c,d 
\in {\cal O},\ c\leq d,$ and (iii) $\hyper {\int_{-\infty}^\infty} f(x) 
\hyper g(x)\, dx \in {\cal O}$ for each $g \in {\cal D}.$\par
\smallskip
\hrule
\vskip 14pt
A function $f$ such that $f^2$ is a limited integral over limited intervals$,$ 
will be said to 
have the {\it limited} (ii) property. 
\vskip 14pt
{\bf Theorem 2.2.} {\sl Suppose that internal $f \colon \hyperreal \to 
\hyperreal,$ has the limited (ii) property$,$ then $f \in T.$} \par
\smallskip Proof. From the hypotheses$,$ $f$ satisfies (i) and (ii) of Definition 2.1.
We need only show that $f$ satisfies (iii). We know that for $g \in {\cal D}$ there 
is some $c \in \real^+$ such that $g(x)= 0$ for all $x \in \real$ such that
$\vert x \vert \geq c.$ Then by Schwarz's inequality (in concise notation)
$$\left(\hyper {\int_{-\infty}^\infty} f \hyper g\right)^2=\left(\hyper {\int_{-
c}^c} f \hyper g\right)^2 \leq $$
$$\left(\hyper {\int_{-c}^c} f^2\right) \left(\hyper {\int_{-c}^c} \hyper g^2\right).$$
But since $g$ is continuous on $[-c,c],\  \left(\int_{-c}^c g^2\right)=r\in 
\real$ 
implies that $\left(\hyper{\int_{-c}^c} \hyper g^2\right)=r \in {\cal O}.$ Since $f$ 
has the limited (ii) property$,$ then 
$$\left(\hyper{\int_{-c}^c} f^2\right) \left(\hyper{\int_{-c}^c} g^2\right)=h\in 
{\cal O}.$$  Consequently$,$ $\hyper{\int_{-c}^c} f\hyper g \in {\cal O}$ and the proof 
is complete.\par
\vskip 14pt
{\bf Example 2.1.} Does $T$ contain nonextended standard functions?
Let $0 \not= \eps \approx 0$ (i.e. $\eps \in \monad 0,$ where $\monad 0$ is 
the set of infinitesimals.) Define the function $f =\{(x,y)\mid (x \in 
\hyperreal)\land (y \in \hyperreal) \land (y = \eps)\}.$ From the internal 
definition principle [3]$,$ $f$ is an internal *-constant function. It is not the extension of a 
standard function since $\eps \notin \real.$ Moreover$,$ $\hyper {\int_c^d} f^2 
= \eps^2(d-c)\in \monad 0 \subset {\cal O}$ for all $c,d,\ (c\leq d)$ such 
that $c,d \in {\cal O}.$ By Theorem 2.2$,$ $f \in T.$  Also consider the 
internal function $f_1 = \{(x,y) \mid (x\in \hyperreal)\land(y \in 
\hyperreal)\land(y = \hyper{\sin}(x + \eps))\}.$ Now $f_1^2$ is *-continuous on 
$\hyperreal$ and$,$ hence$,$ *-integrable on any $[c,d], \ (c\leq d),\ c,d \in 
{\cal O}.$ Further$,$ $-1 \leq \hyper {\sin}^2 (x + \eps) \leq 1$ for all
$x \in \hyperreal.$ Hence$,$ $-1(d-c) \leq \hyper {\int_c^d}\hyper {\sin}^2 
\leq 1(d-c)$ for 
all limited $c\leq d$  implies that $\hyper {\int_c^d}\hyper {\sin}^2 \in {\cal O}$ 
for all limited $c \leq d.$ Again by Theorem 2.2$,$ $f_1 \in T.$
\par
\vskip 14pt
{\bf Theorem 2.3} {\sl The set $\sig (CS) \subset T.$}\par
\smallskip Proof. From the discussion prior to Theorem 2.1.\par
\vskip 14pt   
From   
Theorem 2.1$,$ Corollary 2.1.1$,$ and by Theorem 2.2$,$ $T$ contains many significant 
extended standard functions$,$ among others$,$ that 
will be useful later in this investigation. The next result also indicates 
that one of the conclusions is sufficient for an internal function to be a 
member of $T.$\par
The next 
definition indicates why (iii) in Definition 2.1 is significant. \par
\vskip 14pt
\hrule
\smallskip 
{\bf Definition 2.2 {(The Quasi-Inner Product)}} Let $f \in T$ and $g \in 
{\cal D}.$ Define  $\langle 
f,\hyper g\rangle =  {\hyper{\int_{-\infty}^\infty} f(x) \hyper g(x)\, dx}.$ 
Note that
$\langle \cdot,\cdot \rangle$ is an ordered pair notation since the elements 
come from possibly different sets.\par
\smallskip
\hrule
\vskip 14pt
\noindent {\bf 3. Some Abstract Algebra.}
\vskip 14pt

For most applications of the theory of generalized functions$,$ it may be 
assumed or a function can be appropriately redefined so that the standard 
function being considered is at the least bounded on the closed intervals. 
Call $[c,d]$$,$ where $c \leq d $ and $c,d \in \hyperreal$ a {\it limited *-
closed interval.}
Recall that if $f \in {\cal L}^2(E)$ for measurable
$E$ such that $m(E) \in \real,$ then $f \in {\cal L}(E).$ {[7$,$ p. 220]}   
For the internal 
functions that concern us$,$ *-transfer says that if $f \in \hyper {\cal 
L}^2([c,d]), (c \leq d,\ c,d \in \hyperreal),$ then  
$f \in \hyper {\cal L}[c,d].$
This last statement certainly holds for the limited *-closed intervals. 
Recall the standard result that if $f \in {\cal 
L}([a,b]),$ and $g$ is measurable and bounded on $[a,b],$ then $fg \in {\cal 
L}([a,b]$ [7$,$ p. 219]. Thus by *-transfer for internal function $f$ that is 
*-bounded on a limited *-closed interval $[c,d],$ the function 
$f \in \hyper {{\cal L}^2}([c,d])$ if and 
only if $f \in \hyper{{\cal L}}([c,d]).$ As shown in the next theorem$,$ in some 
special cases$,$ the set $T$ is closed under multiplication of functions.  
\par
\vskip 14pt
{\bf Theorem 3.1} {\sl The following algebraic properties hold. \par
 (a) The set $T$ is an unitary ${\cal O}$- module (left and right) over the 
set of limited numbers 
$\cal O$ and $T$ is a linear space over $\real.$\par
(b) If $f,h \in T$ and $f,h$ are *-bounded in limited *-closed intervals and 
the product $fh$ has the limited (ii) property$,$ then $fh \in T.$\par  
(c) For the set of all continuous functions defined on $\real$$,$ $C(\real),$ 
$\sig C(\real)\subset T$\par
(d) The set $\sig {\cal D}$ is an ideal in ring with unity $\sig C^\infty.$ \par
(e) If $f \in T$ and $\hyper g \in \sig C^\infty,$ then $f\hyper g \in T.$\par
(f) The  the real valued operator  $\langle \cdot,\cdot 
\rangle$ is linear with respect to the field $\real$ in the first and 
second coordinates. The standard part of $\langle\cdot,\cdot \rangle$ is an 
inner product on $\sig {\cal D}.$} \par
Proof.\parm\smallskip     
(a) From the *-transfer of the known 
properties of ${\cal L}^2([c,d]),$ 
$T$ is closed under 
function addition. Since ${\cal O}$ is a ring$,$ if $\lambda \in {\cal O},$ 
then $\lambda f \in T.$ The functions $\Hyper {\b 1} \equiv 1,\ \Hyper {\b 0} 
\equiv 0 \in T.$ 
Hence$,$ $T$ is a unitary ${\cal O}$-module over the ring $\cal 
O.$\parm\smallskip
(b) A standard result says that if $f \in {\cal L}^2([a,b])$ and $g \in 
{\cal L}^2([a,b]),$ then $fg \in {\cal L}([a,b]).$ By *-transfer$,$ this 
statement holds for the limited *-closed intervals. If internal $f$ and 
internal $g$ are 
*-bounded on a limited *-closed interval$,$ then internal $fg$ is *-bounded
on a limited *-closed interval. From our discussion prior to the statement of 
Theorem 3.1$,$ in this case$,$ $fg \in \hyper {{\cal L}}([c,d]).$  But the *-measurable  
internal function 
$fg$ is *-bounded on limited *-closed intervals. Thus  $fg \in \hyper {{\cal 
L}^2}([c,d])$ for the limited *-closed 
intervals. Further$,$ $(fg)^2$ is *-bounded on $[c,d].$ Obviously$,$ 
$fh$ satisfies (i) and$,$ from the hypothesis$,$ (ii) of Definition 2.1. Now by 
Theorem 2.2$,$ it follows that $fh \in T.$ 
\parm\smallskip
(c) $\sig C(\real) \subset \sig (CS)\subset T.$ \parm
(d) The sum and product of members of $C^\infty$ with bounded support have 
bounded support and $\b 0 \in {\cal D}.$ Hence $\sig {\cal D}$ is a subring of $\sig 
C^\infty.$ From
 the bounded support property$,$ if $h \in C^\infty$ and 
$g \in {\cal D},$ then $gh$ has bounded support. Thus $\sig {\cal D}$ is an ideal in
$\sig C^\infty.$\par
(e) Suppose you have a standard  function $f\colon \real \to \real$ such that 
$f\in {\cal L}^2([a,b])$ on standard interval $[a,b].$ 
Since $h^2 \in C^\infty$ is bounded and measurable on $[a,b]$$,$ 
$(fh)^2 \in {\cal L}([a,b])$ by Theorem 22.4s in {[6$,$ p. 127.].}  
Thus$,$ by *-transfer$,$ for $f \in T$ and $\hyper h \in \sig C^\infty,$ $f\hyper 
h$ satisfies (i) and (ii) of 
Definition 2.1. From (d)$,$ if $g \in {\cal D},$ then $hg \in {\cal D}.$ Consequently$,$ 
$\hyper {\int_{-\infty}^\infty} f\hyper h\hyper g \in {\cal O}$ for each $g \in {\cal D}.$ 
\par
(f)  Only the basic algebra for members of $T,$ where it is 
defined$,$ needs to be verified.  Let $f \in T,\  \lambda \in \real.$ 
We know that if $f,h \in T,$ $f+h \in T$ and $\lambda f,\ \lambda h \in T.$ 
If $\hyper g_1,\hyper g_2 \in \sig {\cal D},$ $\hyper g_1 +\hyper g_2 \in \sig {\cal D}$ 
and $\lambda 
\hyper g_1,\ \lambda \hyper g_2 \in \sig {\cal D}.$
For the second coordinate$,$ $\lambda\langle f, \hyper g_1 + \hyper g_2) 
\rangle=
 \lambda\hyper {\int_{-\infty}^\infty} f(\hyper g_1 + \hyper g_2) =
 \hyper{\int_{-\infty}^\infty}\lambda f\hyper g_1 + \lambda f\hyper g_2 =
 \hyper{\int_{-\infty}^\infty}\lambda f\hyper g_1 + 
\hyper{\int_{-\infty}^\infty}\lambda f\hyper g_2=  \hyper{\int_{-
\infty}^\infty}\lambda f\hyper g_1 +\hyper{\int_{-\infty}^\infty}\lambda 
f\hyper g_2=
\lambda \langle f,\hyper g_1 \rangle+ \lambda \langle f,\hyper g_2\rangle.$  
In like 
manner for the first coordinate and$,$ from the above$,$ $\lambda \langle f,g 
\rangle = \langle \lambda f,g \langle = \langle f, \lambda g \rangle.$  
The standard part operator is linear over $\real.$  Hence$,$ the composition of 
$\langle\cdot,\cdot\rangle$ and $\st \cdot$ is linear over $\real.$  Moreover$,$ this composition 
yields a member of $\real.$ Now $\langle \cdot,\cdot \rangle$ is defined on 
all members of $\sig {\cal D}$ independent of order. Further$,$ if 
$g,h \in \sig {\cal D},$ then $ \langle g,h 
\rangle =  \langle h,g \rangle\in \real$ implies that 
$\st {\langle g,h 
\rangle} =  \st {\langle h,g \rangle},$ and $ \langle g,g \rangle \geq 0$ 
implies that $\st {\langle g,g \rangle \geq 0}$.  Since ${\cal D}$ contains only 
continuous functions (with bounded support)$,$ $\langle g,g \rangle = 0$ if and 
only if $g = \b 0.$ Thus on $\sig {\cal D}$ the operator $\st {\langle \cdot,\cdot 
\rangle}$ is an inner product. \par 
\vskip 14pt
\noindent {\bf 4. Functionals on $T.$}\par
\vskip 14pt
\hrule
\smallskip
{\bf Definition 4.1 (The Functional)} Let fixed $f \in T.$ Then$,$ for each 
$ \hyper g 
\in \sig {\cal D},$
define $f[g] = \st {\langle f, \hyper g \rangle}.$ 
Let $T_0= \{f\mid (f \in 
T)\land \forall x((x = \hyper g)\land (\hyper g \in \sig {\cal D}) \to 
(\st {\langle f, x \rangle}=0))\}.$\par
\smallskip
\hrule
\vskip 14pt
From the fact that the standard part operator is linear over $\real,$  it 
follows from Theorem 3.1 (f)$,$ that $f[\cdot]$ is a linear functional on$\sig 
{\cal D}.$
Let ${\cal F} = \{f[\cdot] \mid f \in T\}.$ Obviously there is a surjection
$\Phi \colon T \to {\cal F}.$ Since $T$ is a linear space over $\real$ 
and $\Phi (\b 
0) = {\b 0}[g] = 0$ for each $g \in \sig {\cal D}$  and preserves scalar products 
and ``sums,'' 
$\Phi$ is a vector space homomorphism.\par
\vskip 14pt
{\bf Theorem 4.1.} {\sl The function $f \in T_0$ if and only if $\langle f, \hyper g 
\rangle \approx 0$ for all $g \in {\cal D}.$} \par
\smallskip Proof. This comes from the fact that $\st a = 0$ if and only if $a \approx 0.$\par
\vskip 14pt
{\bf Theorem 4.2} {\sl If $f \in T,$ and $\hyper {\int_c^d} (f(x))^2 \, dx 
\approx 0$ for all limited*-closed intervals$,$ then $f \in T_0$.}\par
\smallskip Proof. This is the same as the proof of Theorem 2.2.\par
\vskip 14pt
Obviously$,$ more than  one $f \in T$ can yield the zero functional and
$T_0 \not= \emptyset$ since $\b 0 \in T_0.$ Now $T_0$ is the kernel for this 
homomorphism$,$ and$,$ as is-known$,$  the quotient linear space $T/T_0$ is 
isomorphic to  ${\cal F}.$ Each element in $T/T_0$ is an equivalence 
class of members of $T.$ Then for $h,g \in T,$ it follows that $f,h \in \alpha \in T/T_0$ if and only if $h - f \in T_0$ if and only if 
$h [g] = f[g] = 0$ for all $\hyper g \in \sig {\cal D}.$  It is this isomorphism 
that allows us to correspond a subset of $T/T_0$ to all of the 
Schwarz generalized functions. \par
\vskip 14pt
{\bf Example 4.1.} Does $T$ contain functions that yield the Dirac property?
Let $$b(t)= \cases{\exp (-1/(1 -t^2)),&$\vert t\vert <1$\cr
                   0,& elsewhere\cr}$$
This is a version of Cauchy's flat function and it is known that $b \in 
{\cal D}.$ We now compress this function. Let $0<\eps\in \monad 0$ (a positive 
infinitesimal). Let $c(t) = b(t/\eps).$ By *-transfer$,$ $c \in \Hyper C^\infty$ 
with support $[-\eps,\eps].$ We can normalize $c$ by letting $k = 
\hyper {\int_{-\infty}^\infty} c(t)\, dt\not= 0,$ and writing $d(t) = 
(1/k)c(t).$ Obviously$,$ $d \in \Hyper C^\infty,$ is nonnegative$,$ and 
$\hyper {\int_{-\infty}^\infty} d(t) \, dt = 1.$ First$,$ we show that 
$d \in T.$ By *-transfer$,$ $d$ is *-bounded *-measurable on each $[c,d],\  c,d 
\in \hyperreal.$ Hence $d \in \hyper {{\cal L}^2}([c,d]),\ c,d \in {\cal O}$ 
(i.e. $d \in \Hyper (CS)$) and satisfies (i) and (ii) of Definition 
2.1. We need only show that for all $\hyper g \in \sig {\cal D},$  $\langle d, \hyper g 
\rangle \approx \hyper g(0)=g(0).$ Recall that the operators $\Hyper \max$ and 
the $\max$ on $\hyperreal$ are the same operator. By *-transfer of the 
standard theorem$,$ 
$$\vert \hyper {\int_{-\infty}^\infty} d \hyper g \vert \leq \sup \{ 
\vert \hyper g(t) \vert\mid t 
\in [-\eps,\eps]\}\cdot \hyper {\int_{-\eps}^\eps} d.$$
Note that $\hyper {\int_{-\eps}^\eps} d = 1.$ The function $\vert g\vert $ is 
continuous at $0.$ Hence$,$ for each $t \in \monad 0,\ 
 \vert \hyper g(t)\vert  \approx g(0).$  From the *-transfer of the extreme value theorem,
there exists some $t_1 \in [-\eps, \eps] \subset \monad 0$ such that 
$\sup \{\vert \hyper g(t)\vert  \mid t 
\in [-\eps,\eps]\} = \vert \hyper g(t_1)\vert$ (i.e. $\sup = \max$). 
Hence$,$ $ \sup \{\vert \hyper g(t)\vert  \mid t 
\in [-\eps,\eps]\} \in \monad {g(0)} \subset {\cal O}.$ 
Thus $d \in T \cap \Hyper (CS).$  In a similar manner$,$ noting that $d$ is 
nonnegative$,$ 
we have we have that 
$g(0) \approx \min \{\hyper g(t) \mid t 
\in [-\eps,\eps]\} \leq \langle d, \hyper g \rangle \leq 
\max \{\hyper g(t) \mid t \in [-\eps,\eps]\} \approx g(0).$ Consequently$,$ the 
functional $d [g] = g(0)$ for all $\hyper g \in \sig {\cal D}.$ But this last 
statement is the  ``shifting'' property of Dirac when viewed as a *-Lebesgue 
integration over $\hyperreal.$ The same method shows that for any positive 
$n \in \nat,\ d^n[g] = g(0).$ \par
\vskip 14pt
\hrule
\smallskip
\hrule
\smallskip
There are infinitely many internal 
functions in 
$\Hyper C^\infty$ that are in $T$ and that determine the Dirac functional. In 
the standard theory$,$ no such standard function exists and ``something'' is 
only 
symbolically introduced relative to the required shifting property. This yields 
what are called ``singular'' generalized functions. From the nonstandard 
viewpoint$,$ at least for the $d[\cdot]$$,$ such a concept of  ``singular''
is no longer meaningful. \par
\smallskip
\hrule
\smallskip
\hrule
\vskip 14pt 
 Since $T/T_0$ is isomorphic to ${\cal F}$$,$ then  
each $f \in T$ such that  $f[g] = g(0)$ for all $g \in {\cal D}$
are in the same member of $T/T_0.$ We call this the {\it Dirac delta}
equivalence class and denote it by $\delta.$ Note that $d^n \in \delta.$\par
\vskip 14pt
{\bf Example 4.2.} The set $\Hyper {\cal D} \not\subset T.$ Consider the function
$b$ of example 4.1. Then $0 <\inf \{d(t)\mid -1/2 \leq t \leq 1/2\}= d(1/2)\in 
\real.$ Further$,$ $0 < \int_{-\infty}^\infty d = r \in \real.$ 
As pointed out$,$ $b \in {\cal D}.$ Let $\Lambda \in \nat_\infty$ (the infinite 
natural numbers).
Then by *-transfer,
$\Lambda b \in \Hyper {\cal D}.$ Now $\inf \{\lambda \Hyper d(t)\mid -1/2 \leq t 
\leq 1/2\}=
\Lambda d(1/2) \in \hyperreal - \real$ and
$$\Lambda d(1/2) \left(\hyper {\int_{-\infty}^\infty} \hyper d \right) \leq 
\hyper \langle \lambda \hyper d, \hyper d \rangle.$$
Thus $\hyper \langle \lambda \hyper d, \hyper d \rangle \notin {\cal O}.$\par
\vskip 14pt\vfil\eject
\hrule
\smallskip
{\bf Definition 4.2 (Pre-generalized Functions)} Each member $\alpha$ in the 
quotient
linear space $T/T_0$ is called a pre-generalized function and each 
member of $T$ is a generalized function. From this 
point on$,$ lower case Greek letters will always denote pre-generalized 
functions. \par
\smallskip
\hrule
\vskip 14pt
One of the reasons$,$ the set $T/T_0$ is called the set of pre-generalized 
functions is 
that a member of $T/T_0$ need not correspond to a Schwarz generalized 
function. But before corresponding pre-generalized functions to Schwarz 
generalized functions$,$ we have the following remarkable result first proved 
by Robinson. The functionals in ${\cal F}$ are specifically generated by the 
$\hyper {\int_{-\infty}^\infty}.$ Does this exhaust the entire collection of all linear 
functionals defined on $\sig {\cal D}$? The following result shows the power of the 
enlargement concept. \par
\vskip 14pt
{\bf Theorem 4.3.} {\sl Let $\Delta$ represent any linear 
functional defined on $\sig {\cal D}.$ Then there exists a *-polynomial $p_\Delta\in \Hyper 
{C^\infty}\cap T$ such that $\Delta = {p_\Delta}[\cdot] \in {\cal 
F}.$}\par
\smallskip Proof. First$,$ let $\Pi$ be the set of all polynomials defined on $\real$ and
$\Delta \colon {\cal D} \to \real.$  
Note that $\Pi \subset C^\infty.$ Consider the binary relation 
$R = \{((g,\Delta(g)),p)\mid (g \in {\cal D})\land(\Delta(g) \in \real)\land
(p\in \Pi)\land (\int_{-\infty}^\infty pg=\Delta(g))\}.$ What is needed 
is to show that $R$ is concurrent on the domain $\real \times {\cal D}.$ \par
Consider  nonempty ${\cal D}$ and a nonempty finite linear independent 
$L=\{g_j\mid (j = 1,\ldots,m)\land (1\leq 
m)\}\subset {\cal D}.$ Then there exists some $c >0$ such that $g_j(x) = 0,\ i\leq 
j\leq m,$  and $c$ can be selected so that each $g_j$ is zero in a 
neighborhood of $-c$ and $c.$ 
 From this we also have that for any
$f \in T,$
$$ f[g_j] = \st {\hyper{\int_{-c}^c} f\hyper g_j} = a_j,\ 1\leq j \leq m.\eqno 
(4.2.1)$$\par
Since each $g_j \in {\cal D} \subset C^\infty,$ let each $g_j$ be represented in 
terms of a series expansion of Legendre polynomials $P_i$ where$,$ using a 
simple transformation of the independent variable$,$  
the $P_i$ have 
been extended to converge on $[-c,c]$ rather than $[-1,1].$  
Hence$,$ $g_j(x) = \sum_{n=0}^\infty a_n^jP_n(x),$ for all $x \in [-c,c]$ and the 
convergence being uniform on any $[-d,d],\ 0<d<c.$ We now use the method of 
the infinite  matrix$,$ a method used by Robinson and Bernstein to solve a 
specific case of the invariant subspace problem. From the linear independent 
assumption$,$ the matrix
$$B=\left(\matrix{a_0^1&a_1^1&a_2^1&\cdots\cr
                \cdot&\cdot&\cdot&\cdots\cr
                \cdot&\cdot&\cdot&\cdots\cr
                \cdot&\cdot&\cdot&\cdots\cr
                a_0^m&a_1^m&a_2^m&\cdots\cr}\right)$$
must be of rank $m.$ Thus there is a finite collect of members $B$ and$,$ hence$,$ 
of subscripts $0\leq j_i\leq\cdots\leq j_m$ such that
$$A = \left(\matrix{a_{j_1}^1&a_{j_2}^1&a_{j_3}^1&a_{j_m}^1\cr
                \cdot&\cdot&\cdot&\cdot\cr
                \cdot&\cdot&\cdot&\cdot\cr
                \cdot&\cdot&\cdot&\cdot\cr
                a_{j_1}^m&a_{j_2}^m&a_{j_3}^m&a_{j_m}^m\cr}\right)$$
and $A$ is nonsingular. Now write
$$A^{-1}\left(\matrix{g_1(x)\cr
                      \cdot\cr
                      \cdot\cr
                      \cdot\cr
                      g_m(x)\cr}\right) =  \left(\matrix{h_1(x)\cr
                      \cdot\cr
                      \cdot\cr
                      \cdot\cr
                      h_m(x)\cr}\right).$$ 
                       
\noindent Thus$,$ we can write $$h_l(x) = P_{j_l}(x) + k_l(x),\ 1\leq l\ \leq m,$$
where the Legendre polynomials in each $k_l(x)$ do not contain any of the
$P_{j_l}(x),\ 1\leq l\leq m.$ \par
Now let
$$A^{-1}\left(\matrix{a_1\cr
                      \cdot\cr
                      \cdot\cr
                      \cdot\cr
                      a_m\cr}\right) =  \left(\matrix{b_1\cr
                      \cdot\cr
                      \cdot\cr
                      \cdot\cr
                      b_m\cr}\right).$$ 
Consider the polynomial $p(x) = c_1P_{j_1} + \cdots c_mP_{j_m}(x).$ We want to 
obtain the proper $c_1,\ldots, c_m$ such that $p[g_j] = a_j,\ 1\leq j 
\leq m.$ First$,$ note that from the orthogonality of the Legendre polynomials 
$\int_{-c}^c pk_l= 0,\ 1\leq l\leq m$ and we know that $\int_{-c}^c P_n^2 =
r_n \not= 0.$ Select $c_l = (1/r_l)b_l,\ 1 \leq l\leq m.$ Substituting
this $p$ with these coefficients into (4.2.1) yields
$$ p[g_j]= \st {\hyper{\int_{-c}^c} \hyper p\hyper g_j} 
= \int_{-c}^c p g_j= a_j,\ 1\leq j \leq m.\eqno (4.2.2)$$\par
\baselineskip=14pt
Since every member of ${\cal D}-L$ (if any) is a linear combination of the members
of $L,$ then it follows from the linearity of $\Delta$ that
for each $g \in {\cal D},\ p[g] = \Delta (g).$ Hence$,$ the relation
$R$ is concurrent on $\real \times {\cal D}$ and$,$ hence on $\sig \real \times 
\sig {\cal D} .$ Thus$,$ there exists some $p_\Delta \in 
\Hyper \Pi \subset \Hyper C^\infty$ such that 
for all $(g,\Delta (g)) \in \real \times \sig {\cal D}, \ ((g,\Delta (g)), 
p_\Delta)\in \hyper R.$ Since $\real\subset {\cal O},$ 
$p_\Delta \in T$ and $\hyper {\int_{-\infty}^\infty} p_{\Delta}\hyper g = \Delta (g).$ 
Consequently$,$ $ {p_\Delta} [g] = \Delta (g)$ and the proof is complete.\par
\vskip 14pt           
Theorem 4.3 is remarkable$,$ since if $\alpha \in T/T_0,$ there is 
$p_\Delta \in \Hyper \Pi$ such that $p_\Delta \in \alpha.$ Furthermore,
$p_\Delta$ is a *-finite sum of *-Legendre polynomials. This means$,$ as 
difficult as it might be to imagine$,$ there is a $p_\Delta \in\delta$ such that 
${p_\Delta}[g] = g(0)$ for each $\hyper g \in\sig {\cal D}.$ The internal function 
$d$ used in Example 
4.1 is NOT a member of $\Hyper \Pi$ by *-transfer of the standard properties 
of the standard function $c(t/a),\ 0<a.$  The function $p_\Delta$ that 
generates every linear functional is also a member of $T_R.$ 
Of course$,$ Theorem 4.3 holds for 
other functions as well that are either simple modifications of the  $P_i$
or are such things as finite sums of trigonometric functions. 
{\bf Note that every linear functional on $\sig {\cal D}$ (i.e. ${\cal D}$) can be generated by the 
Definition 4.1 process by means of *-Riemann integration.} 
Although it is possible to remove 
many functions from each pre-generalized function $\alpha$ by considering the
quotient group formed from the sets $P_0 = T\cap \Hyper C^\infty,\ ({\rm or}\ 
T\cap \Hyper \Pi)$ and $p_1 = T_0\cap \Hyper C^\infty,\ ({\rm or}\ 
T_0\cap \Hyper \Pi)$$,$ it is more significant to have such functions as 
$d \in \delta.$ Later for a necessary simplification process$,$ we will
call a generalized 
function the function in $\Hyper {C^\infty}$ that exists by Theorem 4.3
and not consider the equivalence class at all.\par
\vskip 14pt
\noindent {\bf 5. Schwarz Generalized Functions.}\par
\vskip 14pt
\hrule
\smallskip
{\bf Definition 5.1 (Schwarz Generalized Functions)} A $f[\cdot] \in {\cal F}$ 
is a {\it Schwarz generalized function} if
given any  sequence $\{g_n\} \subset {\cal D}$ such that\par
(i) there exists some $a \in \real$ such that for each $n \in \nat,\ 
g_n(x) = 0$ for all $\vert x \vert >a,$\par
(ii) for each natural number $k \geq 0,\ g_n^{(k)}(x) \to \b 0$ uniformly on
$[-a,a],$ \par
(iii) then   $f[g_n] \to 0.$\par
Let ${\cal D}^\prime$ denote the set of Schwarz generalized functions. \par
\smallskip
\hrule
\vskip 14pt
For every $f[\cdot] \in {\cal D}^\prime,$ there exists unique $\alpha \in T/T_0$
under our isomorphism.  Such an $\alpha$ 
is called
 ${\cal D}^\prime$-pre-generalized function and an $f\in \alpha$ an 
${\cal D}^\prime$-generalized function. 
[The term {\it Schwarz generalized 
function} means the linear functions.]   
As previously pointed out$,$ $\sig (CS) 
\subset T.$ Using only $k = 0,$  properties of Lebesgue integration$,$ 
the Schwarz inequality$,$ if $f \in CS$$,$ then it follows that 
${\hyper f}[\cdot] \in {\cal D}^\prime.$
Using the *-transform of these previous properties,
it follows$,$ using the internal function $d$ defined in Example
4.1$,$ that $d[\cdot] \in {\cal D}^\prime$ although $d$ is not an extended 
standard function. This shows the advantage of selecting specific members of 
a pre-generalized function. \par
\vskip 14pt
{\bf Theorem 5.1.} {\sl Let $f,\ h \in \alpha.$ Suppose that using the 
*-transform of the definition of the derivative that $f$ and $h$  
possess a 
*-derivative $f^\prime$ and $h^\prime$ for each $x \in \hyperreal$ and that 
$f^\prime,\ h^\prime$ are *-continuous for each $x \in \hyperreal.$ Then 
$f^\prime, \ h^\prime \in T$ and there exists a $\beta$ such that 
$f^\prime,h^\prime \in \beta.$}\par
\smallskip Proof. By *-transfer of the continuous case$,$ it follows that 
$f^\prime$ satisfies (i)$,$ (ii) of Definition 4.1.  
By *-integration by parts,
$\st {\langle f^\prime, \hyper g\rangle} = -\st {\langle f, \hyper 
g^\prime\rangle}\in \real$ for each $g \in {\cal D}.$ 
Since $g^\prime \in {\cal D},$ then $f^\prime \in T.$ In like 
manner$,$ $h^\prime \in T.$ Consider $k = f-h.$ We know that $\langle k,\hyper 
g\rangle \approx 0$ for each $g \in {\cal D}.$ Since $k$ satisfies all the 
*-transformed derivative rules$,$  $\st {\langle k^\prime, \hyper g\rangle} = -
\st {\langle k, \hyper g^\prime\rangle}=0$ for each $g \in {\cal D}.$ Thus $k^\prime 
\in T_0$ by our remark after Theorem 4.2. Hence$,$ there exists a unique $\beta 
\in T/T_0$ such that $f^\prime,\ h^\prime\in \beta.$ This completes the 
proof.\par
{\bf Corollary 5.1.1} {\sl For every $f[\cdot] \in {\cal F},$ there exists a 
$f^\prime[\cdot] \in {\cal F}$ such that $f^\prime [g] =- 
f[g^\prime].$}\par           
\vskip 14pt
Since every pre-generalized function $\alpha$ contains an internal $f$ such that there is 
an internal function $f^\prime$ which is the 
*-continuous *-derivative of $f$ on $\hyperreal,$  then to every pre-generalized 
function 
there corresponds a unique $\beta$ such that $f^\prime \in \beta.$ We denote 
such a pre-generalized function by the notation $\alpha^\prime.$ \par
\vskip 14pt
{\bf Theorem 5.2.} {\sl For every $k \in \nat,$  and every $\alpha$$,$ there 
exists a unique $\alpha^{(k)} \in T/T_0.$}\par
\smallskip Proof. Theorem 4.3 shows that for every $\alpha \in T/T_0$ there exists an $
p \in \Hyper C^\infty$ such that $p \in \alpha.$  The result follows by induction 
and Theorem 5.1.\par
\vskip 14pt
\hrule
\smallskip
{\bf Corollary 5.2.1} {\sl For every $0\leq k \in \nat,$ and for every 
$f[\cdot] \in {\cal F},$ there exists a 
$f^{(k)}[\cdot] \in {\cal F}$ such that $f^{(k)} [g] =(-1)^k 
f[g^{(k)}].$}\par
\smallskip
\hrule           
\vskip 14pt
One can now see why included in the definition of the 
Schwarz generalized function 
is the additional part of (ii) for 
each $k >0.$ For$,$ from above$,$ we have  the following important 
result.\par
\vskip 14pt
{\bf Theorem 5.3.} {\sl If $\alpha$ is a ${\cal D}^\prime$-pre-generalized function$,$ then
$\alpha^{(k)}$ is a ${\cal D}^\prime$-pre-generalized function for each $k\in \nat.$}\par
{\bf Corollary 5.3.1.} {\sl If $f[\cdot] \in {\cal D}^\prime,$ then$,$ for each $k \in \nat,$ 
there exists a unique $f^{(k)}[\cdot] \in {\cal D}^\prime$ such that 
$f^{(k)}[g] = (-1)^kf[g^{(k)}].$}\par
\vskip 14pt
Although Corollary 5.3.1 is an important  Schwarz generalized function result$,$ 
the nonstandard theory is more general in that Corollary 5.2.1 holds. \par
\vskip 14pt
\noindent {\bf 6. ``Continuity.''}
\vskip 14pt
Obviously$,$ the definition of a Schwarz generalized function is designed to 
give the linear functional a type of continuity. In nonstandard analysis$,$ 
there are various types of continuity. \par
For two topological spaces $(X,\tau)$ and $(Y,{\cal T}),$   
you always have the 
concept 
of a function $f \colon X \to Y$ as being continuous at $p \in X$ if for each 
$G_1\in 
\tau_1,$ such that $f(p) \in G_1,$  there exists a $G\in \tau$ such that $p 
\in G$ and $f(G) \subset G_1.$ In general$,$ for standard $\hyper p \in \sig X,
$ the 
{\it topological monad} of $\hyper p$ for a given 
topology ${\cal T}$ is $\mu_{\cal T}(\hyper p) = 
\cap \{\Hyper G \mid p\in G \in {\cal T}\}.$  
Then it can be shown that $f$ 
is continuous at $p\in X$ if and only if $\hyper f(\mu_{\tau}(\hyper p)) \subset
\mu_{\cal T}(\hyper f(\hyper p)).$  Note that the reason we need to use 
the standard elements in the form $\hyper p$ is that it is not assumed that 
$X \cup Y$ are atoms within our set-theory. 
Let $\hyperreal^+$ denote the set all positive hyperreal 
numbers. \par
\bigskip
\hrule
\smallskip
{\bf Definition 6.1 (Pseudo-metric Generated Space)} 
Given an internal set $X$ and  $PM_X$ the internal set of all pseudo-metrics 
defined on $X.$ If internal map $\lambda\in PM_X,$ then for each 
$x,\ y,\ z\in X,$ 
(i) $\lambda(x,y)\in \hyperreal^+,$ (ii) $\lambda(x,y) = \lambda(y,x),$ (iii) 
if $x = y,$ then $\lambda(x,y) =0,$ and (iv) 
$\lambda(x,z) \leq \lambda(x,y) + \lambda(y,z).$  Let nonempty $\Lambda \subset
PM_X$. The entity $(X,\Lambda)$ is an 
{\it pseudo-metric generated space.} 
\par
\smallskip
\hrule
\vskip 14pt
Each space $(X,\Lambda)$  satisfies the *-transform of any 
general property for a pseudo-metric. To see this$,$ note that 
there is some 
standard $X_n$ such that $X\in \hyper X_n$ and some standard set $X_p$ 
such that internal 
$PM_X\in \hyper X_p.$ Now$,$ in general$,$ for each $X \in X_n,$ there 
exists a standard set $PM_X.$ Thus there exists a standard set 
${\cal P} = \{PM_X 
\mid (PM_X \in X_p)\land (X \in X_n)\}.$ The internal sets of definition 
6.1 are members of $\hyper X_n$  and the internal $PM_X \in \Hyper {\cal P}.$ 
The defining property for members of internal $PM_X$ is but the *-transform of 
the standard definition. 
Hence$,$ using these sets$,$ any general bounded first-order property about 
standard pseudo-metrics holds$,$ by *-transfer$,$ for members of an internal
$PM_X.$
For example$,$ suppose that internal $\lambda\in \Lambda$  
is determined by an internal 
semi-norm $\Vert \cdot \Vert_\lambda$ defined on an internal space $X$ linear over 
$\hyperreal.$ Then$,$ for $x,y \in  X,$ 
we have that $\vert\  \Vert x \Vert_\lambda - 
\Vert y \Vert_\lambda \ \vert \leq
\Vert x - y \Vert_\lambda \leq \Vert x \Vert_\lambda + \Vert y \Vert_\lambda.$ 
 Of course$,$ our basic examples are the 
standard pseudo-metrics on a standard $X.$ \par
\vskip 14pt
{\bf Example 6.1} Let $SM$ be the set of all internal semi-norms defined
on an internal linear space $X.$ Thus  $\Vert \cdot \Vert \in SM$ 
if and only if for each $
x,\ q \in X,$ (i) $\Vert q \Vert \in \hyperreal^+,$  
(ii) for each $\lambda \in \hyperreal,\ \Vert \lambda q \Vert = \vert 
\lambda \vert\ \Vert q\Vert,$ (iii) $\Vert x + q\Vert \leq \Vert x\Vert 
+\Vert q \Vert.$ Then defining internal $\lambda \colon X \times X \to 
\hyperreal $ by $\lambda(x,q) = \Vert x - q\Vert$ gives $\lambda \in PM_X.$ 
\par
\vskip 14pt\vfil\eject
\hrule
\smallskip
{\bf Definition 6.2. (Monads about $q \in X$)} For  
 $(X,\Lambda)$ and $q \in X$ the {\it monad about q} is
$\mu_\Lambda(q) = \{x\mid (x \in X)\land \forall \lambda((\lambda \in \Lambda) \to 
(\lambda(x,q)\in\monad 0)) \},$  where $\monad 0$  is the set of 
infinitesimals in $\hyperreal.$ \par
Let internal $\lambda \in \Lambda.$ Then the $\lambda$-monad about $q$ is 
$\mu_\lambda(q) = \{x\mid (x \in X)\land (\lambda(x,q)\in\monad 0) \}.$ 
It is an 
important fact that $\mu_\Lambda (q) = \cap \{\mu_\lambda(q)\mid \lambda \in 
\Lambda\}.$           
\par
\smallskip
\hrule
\vskip 14pt
{\bf Example 6.2} Let $\cal S$ be a standard collection of pseudo-metrics 
defined on standard $X.$ Consider the usual topology $\cal T,$  generated by 
the subbase ${\cal B}$ of all open balls determined by all the members of 
$\cal S$ (i.e. $B(p,\lambda,\eps),\ p \in X,\ \lambda \in {\cal S},\ \eps 
\in \real^+$). For a topological space$,$ a topological monad$,$ $\mu_{\cal T}
(\hyper p),$ about standard $p\in X$ is the set $\cap \{\Hyper G \mid p 
\in G \in {\cal T}\}=\cap \{\Hyper G \mid \hyper p 
\in \Hyper G \in \sig {\cal T}\}$ and$,$ in 
general$,$ is equal to $\cap \{\Hyper G \mid p \in G \in {\cal B}\}$ for any 
subbase $\cal B$ for the topology. Thus under Definition 6.2$,$ for $p \in X,$ 
where $\Lambda = \sig {\cal S},$ the monad about $p$$,$ $\mu_\Lambda(\hyper p)$ 
is topological. \par
\vskip 14pt  
\hrule
\smallskip
{\bf Definition 6.3. ($\approx$ and Monads)} 
When an infinitesimal relation  $\approx$ is defined on an internal $X,$ 
this relation is$,$ usually$,$ an equivalence relation and is used to define a 
monad about each $q \in X.$  The monad about $q \in X$ is the 
equivalence class $m_\approx (q) = \{x\mid x \approx q\}$. 
\smallskip
\hrule
\vskip 14pt
In general$,$ the monad defined by 6.3 need not be the same as a monad as 
defined by a topology. But for Definition 6.2$,$ an obvious equivalence relation 
does exist that correspond these monad concepts. \par
\vskip 14pt
\hrule
\smallskip
{\bf Definition 6.4} For the  
space  $(X,\Lambda)$ and any $x,\ y \in X,$  let $x \approx y $ if and only if 
$\forall \lambda((\lambda \in \Lambda) \to 
(\lambda (x,q)\in\monad 0)).$ Also$,$ for each $\lambda \in \Lambda,$ 
$x {\buildrel\lambda \over \approx} y$ if and only if 
$\lambda (x,q)\in\monad 0.$\par
\smallskip
\hrule
\bigskip
It is immediate that the $\approx$ [resp. ${\buildrel\lambda \over \approx}$] 
of Definition 6.4 is  an equivalence 
relation on $X \times X,$ and that $x \approx y$ [resp. $x {\buildrel\lambda 
\over \approx} y$] if and only if $x \in \mu_\Lambda (y)$ [resp. $x \in \mu_\lambda 
(q)]$ if and only if $x \in m_\approx(y)$ [resp. $ x \in 
m_{{\buildrel\lambda \over \approx}}(y)$.] 
\bigskip
 \hrule
\smallskip
{\bf Definition 6.5. (S-continuous)} 
For spaces  $(X,\Lambda)$ and 
$(Y, \Pi)$
an internal $f\colon X 
\to  Y$ 
is {\it S-continuous} at 
 $q \in X$ if $f(\mu_\Lambda(q)) \subset 
\mu_\Pi(f(q)).$ For  pseudo-metric space $(X,\lambda)$ and $(Y,\pi),$
S-continuity is defined for the spaces $(X,\{\lambda\})$ and $(X,\{\pi\}).$\par
\smallskip
\hrule
\bigskip
For an pseudo-metric$,$ $\lambda$ defined on internal $X,$ you 
define$,$ for each $\eps \in \hyperreal^+,$ and for each $q \in X,$ (in the 
usual way) the 
ball about $q$ as $B(q,\lambda,\eps)=\{x\mid (\lambda(x,q) < \eps \}.$ 
Another type of continuity$,$ that is usually restricted to standard
spaces$,$ is *-continuity. For a standard  
gauge space (i.e. the topological space generated by the set of pseudo-metrics
$\Lambda$)$,$ then the topology is generated by taking as a bases the finite 
intersection of standard balls.  Since any finite set of real numbers has a 
minimum$,$ using the neighborhood bases about a standard point$,$ we have that 
for standard $(X,\Lambda),\ (Y,\Pi)$  a standard functions $f\colon X \to Y$ 
is continuous at $p \in X$ if and only if for each $\eps \in \real^+$ and each 
$\pi \in \Pi$ there exists a $\delta \in \real^+$ and a finite set of
pseudo-metrics $\lambda_i \in \Lambda,\ 1 \leq i\leq n$ such that whenever
$x \in X$ and $\lambda_i(x,p) < \delta$ for each $i, \ 1\leq i \leq n$ then 
$\pi(f(x),f(p)) < \eps.$ 
The *-transfer 
of this statement is used to define another type of continuity.
\bigskip  
\hrule
\smallskip
{\bf Definition 6.6 (*-continuity)} Consider  pseudo-metric generated 
spaces $(X,\Lambda),\  (Y,\Pi).$  
An internal map $f\colon ( X, \Lambda) \to ( Y,\Pi)$ is 
*-{\it continuous} at $q \in X$ if for each $\eps \in \hyperreal^+$ and 
each $\pi \in \Pi$ there exists a $\delta \in \hyperreal^+$ and a *-finite set 
$\{\lambda_i\mid 1 \leq i \leq \omega\}\subset \Pi$ such that whenever $x \in X$ and $
\lambda_i (x ,q) < \delta$ for each $i$ such that $1 \leq i \leq \omega$ then
$\pi (f(x),f(q))<\eps.$ 
For internal metric spaces$,$ $(X,\lambda)$ and $(Y,\pi),$ 
*-continuity is defined for $(X,\{\lambda\})$ and $(Y,\{\pi\}).$ 
\par
\smallskip
\hrule
\vskip 14pt
Notice that *-continuity is not defined solely in terms of the  
standard points  
in $X.$  Also the specific $\eps$ and $\delta$ required for *-continuity are 
members of $\hyperreal^+,$ not just the standard positive reals.  
For standard $(X,\Lambda),$ $(Y,\Pi)$$,$  
$f \colon (X, \Lambda)    
 \to (Y,\Pi)$ is continuous at $p \in X$ 
if and only if 
$\hyper f \colon (\hyper X, \hyper \Lambda)    
 \to (\Hyper Y,\Hyper \Pi)$  is *-continuous at $\hyper p$.\par  
\medskip
{\bf Example 6.4.} By *-transfer$,$ for any infinite $\Lambda \in \nat_\infty$,
the function $f(x) = \Hyper \sin (\Lambda x)$ is *-continuous on 
$\hyperreal$ with respect to 
the standard norm $\vert \cdot \vert.$  
\par
\medskip 
The *-continuous functions have an important place in theoretical physics 
[4] [See example 6.6 below]. But since they have all of 
the properties of the continuous
functions$,$ they are not considered 
``interesting'' to some members of the mathematics community.  \par
\vskip 14pt
{\bf Example 6.5.} Although the *-continuous function of example 6.4 is 
*-continuous at $x = 0,$ it is 
S-discontinuous at $x = 0.$  Let infinitesimal $\eps = (1/\Lambda)(\pi/2).$ 
Then for S-continuity we must have that $\Hyper \sin (\Lambda\cdot 0)
=\Hyper \sin (0)= 0 \approx \Hyper \sin (\Lambda ((1/\Lambda)(\pi/2)) =
1$; a contradiction. \par
\vskip 14pt
{\bf Example 6.6.} The function $d$ of example 4.1 that generates the Dirac functional is
*-continuous at $x = 0$ but is S-discontinuous. Since the defining 
statement for $d$ is the *-transform of the collection of all functions
$c$ formed by letting $t \to t/a,$ where $a$ is any nonzero real number$,$ it 
follows that $d$ is *-continuous at all $x \in \hyperreal$ (i.e. $d$ is a 
member of a set of *-continuous functions with this property.) However,
$[-\eps, \eps] \subset d(\monad 0)$ and $d[-\eps,\eps] = \Hyper [0,e^{-1}]$ imply 
that $d$ is not S-continuous at $x = 0.$\par
\vskip 14pt   
 The  ``S'' in S-continuous means 
``standardly'' in the sense that the approximating numbers $\eps, \ \delta$ 
are standard 
numbers. One example shows what S-continuity is trying to accomplish.\par 
\vskip 14pt
{\bf Example 6.7.} Let $0<\eps\in \monad 0$ and $a \in \real.$ 
Define on $\hyperreal$ $$f(x) = \cases{\eps +a,&$  x< 0$\cr   
                 a,&$x \geq b$\cr}$$
This is$,$ by *-transfer$,$ an internal function that is *-continuous for all
nonzero $x \in \hyperreal$ and is *-discontinuous at $x = b.$  But$,$ $f$ 
is S-continuous at $x = b.$ For take any positive $ c \in \real,$ then no matter 
what $x \in \hyperreal$ you select$,$ $\vert f(x)-a \vert < c.$ { \bf That is the 
*-discontinuity is  so ``small'' that it is not ``visible'' in the standard
world.} \par
\vskip 14pt
For the topological spaces used in the theory of generalized functions$,$ does 
the concept of S-continuous correspond to the concept described in Example 
6.7? In order to examine the relation between S-continuity and *-continuity for 
generalized functions$,$ a slight diversion is necessary\par

 For a given standard $X,$ suppose  ${\cal V}$ is a 
collection of subsets of $X$ such that 
$\emptyset \notin {\cal V}$ and ${\cal V}$ has the  finite intersection 
property (i.e. the intersection of finitely many members is not the empty set). 
Then the collection $\cal V$ is a {\it filter subbase} on $X.$ Further$,$ if
there exists some $q \in \hyper X$ such that $q \in \Hyper G$ for each 
$G\in {\cal V},$ then
${\cal V}$ 
is$,$ obviously$,$ a filter subbase which is termed the  
{\it local filter subbase} at $q$ and denoted by ${\cal V}_q.$   
\par
\vskip 14pt
\hrule
\smallskip
{\bf Definition 6.7 (Monads of a Filter Subbase)} For standard $X,$ let 
$\cal V$ be a filter subbase (either local or otherwise) on $X.$   
 Then let
$\mu_{\cal V} = \cap \{\Hyper V\mid V\in {\cal V}\}.$  If $\cal V$ is local 
at $q \in \hyper X,$ then since $q \in \mu_{\cal V}$ the  monad is written as $
\mu_{\cal V}(q)$.\par
\smallskip
\hrule
\vskip 14pt
{\bf Example 6.8} For any standard pseudo-metric space $(X,{\Lambda}),$  
and $p \in X,$ 
consider the set$,$ ${\cal B} = 
\{B(p,\lambda, \eps)\mid (\lambda \in {\Lambda})\land(\eps \in  
\real^+)\},$ of all balls about $p.$  Then $\cal B$ is a local filter 
subbase at  $p.$  For any filter subbase  $\cal B$$,$ let $\langle {\cal B} 
\rangle$ be the set obtained by taking finite intersections of members of 
${\cal B}.$  Obviously$,$ ${\cal B} \subset \langle {\cal B} 
\rangle$ and $\mu_B = \mu_{\langle {\cal B} 
\rangle}.$ \par

\vskip 14pt
{\bf Theorem 6.2.} {\sl Let $X$ be standard set and  $\cal V$ any standard 
filter subbase on $X.$ Then $\mu_{\cal V} \not= \emptyset.$}\par
\smallskip Proof. We know that there is some $X_n,\ n\leq 1,$ such that ${\cal V}\in X_n,$ and if 
$x \in {\cal V},$ then $x \in X_n.$ Further$,$ if $y \in x,$ then $y \in X_0 
\cup X_{n-1}.$
Consider the bounded binary relation 
$$\Phi (x,y) =\{(x,y) \mid (y \in X_0 \cup X_{n-1})\land(x \in X_n)\land
 (y \in x)\land (x \in {\cal V})\}.$$
The domain of $\Phi$ is ${\cal V}.$ 
Let $\{(a_1,b_1),\ldots (a_n,b_n)\} \subset \Phi.$  Since ${\cal V}$ has the 
finite intersection property$,$ there exists some $b\in X_0\cup X_{n-1}$ such 
that $b \in a_1 \cap \cdots \cap 
a_n.$ Thus $\{(a_1,b),\ldots (a_n,b)\} \subset \Phi.$ Hence since we are in an 
enlargement$,$ there exists $a \in \hyper X_0 \cup \hyper X_{n-1}$ such that 
$a \in \Hyper V$ for each $V \in {\cal V}.$ Consequently$,$ $\mu_{\cal V} \not= 
\emptyset$ and the proof is complete.\par
\vskip 14pt
Note that if $\cal V$ is a filter subbase$,$ that the set obtained by taking 
finite intersections of members of $\cal V$ does not contain the empty set and 
is closed under finite intersection.\par
\vskip 14pt
{\bf Theorem 6.3.} {\sl Let $X$ be standard set and  $\cal V$ any standard 
collection of subsets of $X,$ which does not contain the empty set and which 
is closed under finite intersection. If internal $A \subset\mu_{\cal V},$ then 
there  
exists some *-finite (and$,$ hence$,$ internal) 
$\Omega \subset \Hyper {{\cal V}}$ such that 
$\sig {\cal V} \subset \Omega$ and $A \subset A_0=\cap \{B\mid B 
\in \Omega \} \subset \mu_{\cal V},$ and $A_0 \in \Hyper {{\cal V}}.$ }\par
\smallskip Proof. First$,$ from Theorem 6.2 $\mu_{\cal V} \not= \emptyset.$ Let 
$B=\{V\mid (V\in \Hyper {{\cal V}})\land (A \subset V)\}.$ Then $B$ 
is an internal subset of $\Hyper {{\cal V}}.$ Extending Theorem 4.3.4 [3]$,$ we 
know that there exists a *-finite set $B_0$ such that 
$\sig {\cal V} \subset B_0 \subset \Hyper {\cal V}.$ Let $\Omega = B \cap 
B_0.$ Since every internal subset of a *-finite set is *-finite$,$ $\Omega$ is 
*-finite. Further$,$ $\sig {\cal V} \subset \Omega$ and $A \subset E$ for each 
$E\in \Omega.$ By *-transfer$,$ $\Hyper {\cal V}$ is closed under 
*-finite  intersection. Hence $A_0=\cap \{B\mid B 
\in \Omega \}\in \Hyper {\cal V},$ $A \subset A_0$ and $A_0 \subset \mu_{\cal V}$ and 
the proof is complete.\par
{\bf Corollary 6.3.1.} {\sl  Let $X$ be standard set and  $\cal V$ any standard 
collection of subsets of $X,$ which does not contain the empty set and which 
is closed under finite intersection. Then $\mu_{\cal V} = 
\cup \{G\mid (G \in \Hyper {\cal V})\land (G \subset \mu_{\cal V})\}.$} \par
\smallskip Proof. For every $q \in \mu_{\cal V,}$ there is some $G \in \Hyper {\cal 
V}$ such that $q\in G \subset \mu_{\cal V}.$  \par
\vskip 14pt
A nonempty collection $\cal B$ of subsets if $X$ is a filter base$,$ 
if $\emptyset \notin
{\cal B}$ and if $A,\ B \in {\cal B},$  then there exists some $C \in {\cal 
B}$ such that $C \subset A \cap B.$ A filter base is a filter subbase. \par
\vskip 14pt
{\bf Theorem 6.4.} {\sl For a standard set $X,$ let $\cal B$ be a 
standard filter base defined on $X.$  If internal $A \subset \mu_{\cal B},$ then 
there  
exists some *-finite (and$,$ hence$,$ internal) 
$\Omega \subset \Hyper {\langle {\cal B}\rangle}$ such that 
$\sig \langle{\cal B\rangle} \subset \Omega$ and $A \subset A_0=\cap \{B\mid B 
\in \Omega \} \subset \mu_{\cal B},$ and $A_0 \in \Hyper {\langle{\cal 
B}\rangle}.$ }\par
{\bf Corollary 6.4.1.} {\sl For a standard set $X,$ let $\cal B$ be a 
standard filter base defined on $X.$  
Then $\mu_{\cal B} = \cup \{G\mid (G \in \hyper {\langle{\cal 
B}\rangle})\land (G \subset \mu_{{\cal B}})\}.$} \par
\vskip 14pt
{\bf Theorem 6.5.} {\sl Given an internal set $A$ and a standard filter subbase
$\cal F$ such that $A \cap \Hyper F \not=\emptyset$ for each $F \in {\cal F}.$ 
Then $A \cap \mu_{\cal F} \not=\emptyset.$}\par
\smallskip Proof. 
Consider the internal binary relation on ${\cal B} = \langle {\cal F} 
\rangle.$ 
$$\Phi =\{(a,b)\mid (b\in A)\land (b \in a)\land (a \in \Hyper {\cal B})\}.$$ 
Suppose that $\{(a_1,b_1),\ldots,(a_n,b_n)\} \subset \Phi$ and $a_i \in \sig 
{\cal B}.$ Since $\cal B$ is closed under finite intersection$,$ we have  
nonempty 
$a \in {\cal B},$ 
where $\Hyper a = \Hyper b_1 \cap \cdots \Hyper b_n.$ Thus there is some 
internal $b^\prime \in \Hyper a$ such that 
$(a_1,b^\prime),\ldots,(a_n,b^\prime) \subset \Phi.$ Since the cardinality of 
$\cal B$  less than $\kappa,$ then $\Phi$ is$,$ at least$,$ concurrent on $\sig {\cal B}.$ Thus there exists some $q \in A \cap \Hyper B$ for each $B 
\in {\cal B}.$ This implies that $A \cap \mu_{\cal B}=A \cap \mu_{\cal F}
 \not=\emptyset$ and the 
proof is complete.\par  
\vskip 14pt
{\bf Theorem 6.6.} {\sl Let $\cal B$ be a standard filter base and suppose 
that an 
internal set $\Lambda \subset \Hyper {\cal B}$  has the property that 
$\sig {\cal B} \subset \Lambda.$ Then there exists some internal $A \in \Lambda$ such 
that $A \subset \mu_{\cal B}.$ }\par
\smallskip Proof. Consider the internal (bounded) binary relation $$
\Phi = \{(b,a) \mid ((b \in \Hyper {\cal B}) \land( a \in \Lambda) \land(a 
\subset b)\}.$$
Let $\{(b_1,a_1)\ldots,(b_n,a_n)\} \subset \Phi$ and $a_i \sig {\cal B}.$ Then 
there is some $b^\prime \in \sig {\cal B}$ such that $b^\prime \subset b_1 
\cap \cdots \cap b_n.$ But $b^\prime \in \Lambda.$ Hence$,$ 
$\{(b_1,b^\prime)\ldots,(b_n,b^\prime)\} \subset \Phi.$  The 
$\kappa$-saturation$,$ there exists some internal 
$A\in \Lambda$ such that $A \subset \Hyper B$ for each $B \in {\cal B}.$ 
Hence$,$ $A \subset \mu_{\cal B}$ and the proof is complete.\par
\vskip 14pt
{\bf Theorem 6.7} {\sl Let $\cal B$ be a standard filter base and suppose 
that an 
internal set $\Lambda \subset \Hyper {\cal B}$  has the property that  
if $G \in \Hyper {\cal B}$ and $G \subset \mu_{\cal B},$ then 
$G \in \Lambda.$  Then there exists some $B \in {\cal B}$  such that
$\Hyper B \in \Lambda.$}\par
\smallskip Proof. Assume the hypothesis and that there is no $B \in {\cal B}$ such 
that $\hyper B \in \Lambda.$  Then the internal set $\Hyper {\cal B} - \Lambda   
\subset \Hyper {\cal B}$ satisfies the hypothesis for the  ``$\Lambda$'' 
of Theorem 6.6. Thus there exists some $A \in \Hyper {\cal B} - \Lambda$   
such that $A \subset \mu_{\cal B}.$ But from the hypothesis of this theorem,
such an $A \in \Lambda.$ This contradiction complete the proof. \par
\vskip 14pt
In order to obtain a significant result that characterizes S-continuity$,$ we 
need the following additional fact.\par
\vskip 14pt
{\bf Theorem 6.8} {\sl Consider standard $X.$ Let internal ${\cal B} \subset
\Hyper F(\hyper X),$ where $\Hyper F(\hyper X)$ is the set of all 
*-finite subsets of $\hyper X.$ Suppose 
that whenever $E \in   \Hyper F(\hyper X)$ and $\sig X \subset E \subset 
\hyper X$ then $E \in {\cal B}.$  Then there exists $F \in F(X)$ such that 
$\Hyper F\in {\cal B}.$}\par
\smallskip Proof. To establish this$,$ let ${\cal F} = \{F\mid X-F \in F(X)\}.$ Then $\cal 
F$ is a filter on $X.$  This theorem is but 
an equivalent statement of theorem 6.6 in terms of  the filter $\cal F.$ 
Establishing this fact completes the proof. \par
\vskip 14pt
We are now able to properly characterize the concept of S-continuity relative 
to standard pseudo-metric generated spaces.\par
\vskip 14pt
{\bf Theorem 6.9.} {\sl For standard $X,\ Y,$ consider the pseudo-metric
generated spaces $(\hyper X,\sig \Lambda),\  (\hyper Y,\sig \Pi).$ An internal 
$f\colon \hyper X \to \Hyper Y$ is S-continuous at $q \in \hyper X$ if and 
only if for each $\eps \in \real^+$ and each $\hyper \pi \in  \sig \Pi$ there 
exists a finite set $\hyper \lambda_i,\ 1 \leq i \leq n$ and positive $\delta 
\in \real^+$ such that whenever $x \in\hyper X$ and $\hyper \lambda_i 
(x,q)<\delta$ for each $i,\ 1 \leq i\leq n,$ then 
$\hyper \pi (f(x),f(q)) < \eps.$}\par
\smallskip Proof. $\Rightarrow$ First$,$ suppose that $f$ is S-continuous at $q \in \hyper X,$ $\hyper 
\pi \in \sig \Pi$ and let $\eps \in \real^+.$   
For any internal binary relation $A,$ let $D(A)$ denote the internal domain 
and 
$R(A)$ the internal range.  
Consider the internal set
$$T(\eps) =\{K \mid (\emptyset \notin D(K)) \land
(\emptyset \notin R(K))\land (K \in \Hyper F(\hyper \Lambda \times 
\hyperreal^+)\land 
(\forall \lambda\forall \delta\forall x$$
$$((\lambda \in D(K))\land(\delta \in R(K)) \land  
(x \in \hyper X)\land(\lambda(x,q) < \delta) \to $$ $$(\hyper \pi (f(x),f(q))<\eps 
)))\}.$$
By $\kappa$-saturation$,$ we know that there exists a *-finite $K_0 \subset 
\hyper \Lambda \times \hyperreal^+$  
such that 
$\sig (\Lambda \times \real^+)\subset K_0.$  Thus 
$\emptyset \not= T(\eps).$ Further$,$ 
suppose that $K_1 \in \Hyper F(\hyper \Lambda \times \hyperreal^+)$ and $\sig 
\Lambda \times \real^+ \subset K_1.$ Then  
 $D(K_1) =G_1 \in \Hyper F(\hyper \Lambda),\ $  
$\sig \Lambda \subset G_1$ and $R(K) = H_1 \in \Hyper F(\real^+),\ $  
$\real^+ \subset H_1.$ Now $x \approx q$ implies that 
$\hyper \lambda(x,q) < \delta$ for each $\hyper \lambda \in \sig \Lambda$ and 
each $\delta \in \real^+.$ But$,$ S-continuity implies that $\hyper \pi(x,q) 
<\eps.$  Hence $K_1 \in T(\eps).$ Since $T(\eps)$ is internal$,$ then Theorem 
6.8 implies that there exists standard $K^\prime$ such that $K^\prime \in 
F(\Lambda \times \real^+)$ and such that $\Hyper K^\prime \in T(\eps).$
Consequently$,$ there is positive $n\in \nat$ and $\lambda_i\in D(K^\prime) 
\subset \Lambda$ when $1 
\leq i \leq 
n$ and a positive $m$ and $\delta_j\in R(K^\prime) \subset  \real^+$ when 
$1 \leq j\leq m$ and for each $x \in X,$ if 
for each $i$ and for each $j$ $\hyper \lambda_i(x,q)<\delta_j,$  then
$\hyper \pi (f(x),f(q) < \eps.$ Now simply consider $\delta = \min 
\{\delta_1,\ldots,\delta_m \}$ and $\Rightarrow$ holds. \par
$\Leftarrow$  
Suppose that $f$ is not S-continuous at $q.$ Then there exists 
some $x_0\in \mu_{\sig\Lambda}(q)=\cap\{\mu_{\hyper \pi}(q)\mid \hyper \pi \in 
\sig \Pi\},$ such that $f(x_0) \notin \mu_{\hyper \Pi}(f(q)).$ 
Thus there is some $\pi \in \Pi$ such that $f(x_0) \notin \mu_\pi(f(q)).$ Let 
standard $\eps = \min\{1, \st {\pi( f(x_0) ,f(q))}/2\}$ if  
$\pi( f(x_0),f(q)) \in {\cal O},$ in which case $\eps >0$ since
$\st {\pi( f(x_0) ,f(q)}/2 \not= 0.$ Otherwise take $\eps = 1.$
Now for all
standard $\delta >0,$  we have that for each $\hyper \lambda \in \sig \Lambda,\ 
\lambda( x_0 ,q)< \delta,$ but $\hyper \pi( f(x_0) ,f(q)) \geq \eps.$ 
The proof is complete.\par 
{\bf Corollary 6.9.1.} {\sl For standard $X,\ Y,$ consider the pseudo-metric
generated spaces $(\hyper X,\sig \Lambda),\  (\hyper Y,\sig \Pi).$ An extended 
standard map
$\hyper f\colon \hyper X \to \Hyper Y$ is S-continuous at $\hyper p \in \sig X$
 if and 
only if $f\colon (X,\Lambda) \to (Y,\Pi)$ is at continuous at $p \in X.$}\par
\vskip 14pt
Notice that the $\delta$s and $\eps$s that appear on Theorem 6.9 are standard 
real numbers. This theorem shows that close relationship between the concept 
of S-continuity and the concept of continuity. For the only difference within
our $\kappa$ saturated model between these two concepts when viewed from the 
external nonstandard physical world is that one  hand the $x$s are members 
of $\hyper X$ while on the other they elements of $X.$ \par
One of the major facts about S-continuous functions is found in Theorem 1.1 
in [4$,$ p. 805.] As pointed out$,$ for this theorem$,$ the topological space $X$ 
need not be compact and the first two parts of the theorem hold.  
In this theorem$,$ the term microcontinuous is  equivalent to
S-continuous. For simplicity of notation$,$ the sets $X, \ Y$ are  
considered as a subset of the ground set that is used to generate our 
ultraproduct  structure. \par
\vskip 14pt
{\bf Theorem 6.10.} {\sl Suppose that you have the topological spaces  $X,\ Y,$
$Y$ regular Hausdorff$,$ where topological monads are defined at standard points,
and an internal $f \colon \hyper X \to \Hyper Y.$ At $p \in X,$ let $f$ be 
S-continuous 
and suppose that there exists some $ r \in Y$ such that $f(\hyper p) \approx r.$  
Then any function $F\colon X \to Y$ such that $F(p) = \st {f(\hyper p)}$ 
is continuous at $p$ in the topological sense.  Further$,$ if $q \approx \hyper 
p,$ then $f(q) \approx  \Hyper F(\hyper p).$} 
[Note. The standard part operator as defined 
in this theorem$,$ brings points all the way back to the original standard set
$Y.$ It is not  significant that this operator can be considered as defined
on  $Z \subset \hyper Y$ and $\st {Z} \subset \sig Y.$] \par
\vskip 14pt
Theorem 6.10 shows$,$ in all generality$,$ how S-continuity in general topological 
spaces leads directly to a standard function because the internal function
ignores infinitesimal discontinuities. If the space $Y$ is not be Hausdorff,
then$,$ by the Axiom of Choice$,$  continuous function(s) can still be 
constructed. \par
\vskip 14pt
\hrule
\smallskip
{\bf Definition 6.8 (S-convergence)} Let internal sequence $s\colon \hypernat \to
(X,\Lambda),$  where  $(X,\Lambda)$ is a  pseudo-metric generated space.
Then $s$ S-converges to $q \in X$ if $s_\omega \in \mu_\Lambda (q)$ for each 
infinite $ \omega \in \nat_\infty.$\par
\smallskip
\hrule
\vskip 14pt \vfil\eject
{\bf Theorem 6.11.} {\sl If$,$ for  pseudo-metric generated spaces,
 internal $f\colon (X,\Lambda) \to (Y,\Pi)$ is S-continuous at 
$q \in X$ and the internal sequence  $s\colon \hypernat \to
(X,\Pi)$ S-converges to $q,$ then the internal sequence $f(s)$ S-converges 
to $f(q).$}\par
\smallskip Proof. Suppose that the internal sequence $s$ S-converges to $q \in X.$ 
Let $\omega \in \nat_\infty.$  Then $s_\omega \in \mu_\Lambda (q).$ 
Thus$,$ $f(s_\omega) = 
(fs)_\omega \in \mu_\Pi (f(q))$ and the proof is complete.\par
\vskip 14pt
The next theorem is similar to Theorem 6.9$,$ relates S-convergence to standard 
approximations as well as to standard sequences.\par
\vskip 14pt
{\bf Theorem 6.12.} {\sl For a pseudo-metric generated space $(X,\Lambda),$  
an internal sequence $s\colon \hypernat \to X,$ S-converges to $q \in X$ 
if and only if for each positive $\eps \in \real$  
each $\lambda\in \Lambda$ there exists some $M \in \nat$ 
such that for each $m > M$ (in $\hypernat$)$,$ it follows that
$\lambda(s_m,q) <\eps.$}\par
\smallskip Proof. ($\Rightarrow$) Suppose the internal sequence $s$ is S-convergent to 
$q \in X.$ Let 
$0<\eps \in\real.$ Consider the internal set
$$m(\eps) =\{m \mid (m \in \hypernat) \land(\forall n((n \in \hypernat)\land(n> m) 
\to \lambda(s_n,q) <\eps)))\}.$$
From the definition of S-convergence $\nat_\infty \subset m(\eps).$ However,
the set $m(\eps)$ has an internal range. Hence$,$ there exists some standard
$m \in \nat$ such that $m \in m(\eps).$  
Consequently$,$ 
the conclusion follows.\par
($\Leftarrow$) Suppose that $s$ is not S-convergent to $q.$ Then there exists 
some $s_\omega, \  \omega \in \nat_\infty$ such that $s_\omega \notin 
\mu_\Lambda(q).$
Hence$,$ there is some $\lambda \in \Lambda$ such that $s_\omega \notin 
\mu_\lambda(q).$ 
Let standard $\eps = \min\{1, \st {\lambda(s_\omega ,q)}/2\}$ if  
$\lambda(s_\omega,q) \in {\cal O},$ in which case $\eps >0$ since
$\st {\lambda(s_\omega ,q)/2} \not= 0.$ Otherwise take $\eps = 1.$
Thus there exists some $\omega \in \nat_\infty$ such that 
 and $\lambda( s_\omega ,q) \geq \eps.$ The proof is complete.\par
{\bf Corollary 6.12.1.} {\sl For a standard pseudo-metric space $(X, \rho_X),$ 
a standard sequence $s \colon \nat \to (X,\rho_X)$ 
is convergent to $p\in X$  if and only if
$\hyper s$ is S-convergent to $p.$}\par
\vskip 14pt
Now we need to define the limited points for internal pseudo-metrics.  For
an internal pseudo-metric $\rho_X$$,$ 
you have a set of limited points per $p \in X.$  
\vskip 14pt
\hrule
\smallskip
{\bf Definition 6.9 (Limited Points for Pseudo-metric Generated 
Spaces)} Let $(X,\Lambda)$ be a pseudo-metric generated space. Let $q \in 
X$ Then ${\cal O}_\Lambda (q) = \{x \mid (x \in X)\land\forall \lambda((\lambda\in 
\Lambda) \to  
\lambda(x,q) \in {\cal O}
)\}.$ \par
\smallskip
\hrule
\vskip 14pt \vfil\eject
{\bf Theorem 6.13.} {\sl Let $(\hyper X,\Lambda)$ be a pseudo-metric generated space. 
Suppose that each $\lambda \in \Lambda$ is determined by an internal 
pseudo-norm $\Vert \cdot \Vert_\lambda$ and that $\hyper X$ is an internal linear 
space over $\hyperreal.$   
Suppose that for each $p \in X$ and each $\Vert \cdot\Vert_\lambda,
\ \Vert p \Vert_\lambda \in 
{\cal O}.$ Then for each $p \in X,\ {\cal O}_\Lambda(p)= {\cal 
O}_\Lambda(\b 0)= \{x \mid (x \in \hyper X)\land\forall \lambda((\lambda \in 
\Lambda) \to (\Vert x \Vert_\lambda \in {\cal 
O}\}= {\cal O}_X.$}\par 
\smallskip Proof. Let $x \in {\cal O}_\Lambda(p)$ for 
$p \in X.$ Then for any $\lambda \in \Lambda$ it follows that 
$\vert\ \Vert x\Vert_\lambda -  \Vert p \Vert_\lambda\ \vert 
\in {\cal O}.$ Hence$,$ $\Vert x\Vert_\lambda \in {\cal O}$ since $\Vert p 
\Vert_\lambda \in {\cal O}.$\par
Conversely$,$ let $x \in {\cal O}_\Lambda (\b 0)$ and $p \in X.$ Then 
for each $\lambda \in \Lambda,\ \Vert x \Vert_\lambda \in {\cal O}.$ But 
$ \Vert p \Vert_\lambda \in {\cal O}$ implies,
since $\Vert x - p \Vert_\lambda \leq \Vert x\Vert_\lambda + \Vert p 
\Vert_\lambda 
\in {\cal O},$ that $x \in {\cal O}_\Lambda (p).$  The proof is complete.\par
{\bf Corollary 6.6.1.} {\sl Let the standard pseudo-metric $\rho_X$ be defined 
on $X\times X$ 
and be generated by 
a standard semi-norm
$\Vert \cdot \Vert_X.$ 
If $r,\ p \in X$$,$ then ${\cal O}_\rho(p) = {\cal O}_\rho(r) = {\cal 
O}_X.$}\par
\vskip 14pt
By Robinson's Theorem 4.3$,$ the set $T$ can be replaced by the set $\Hyper 
C^\infty \cap T.$ Hence$,$ the usual practice has been to consider defining
sets of internal semi-norms on the set $\Hyper C^\infty$ and consider the 
restriction such a set of internal semi-norms to the test space $D.$ This 
collection can be composed of the nonstandard extensions of the customary 
set of standard semi-norms so that they correspond to the concept of 
Schwartz generalized functions and other standard types of generalized 
functions. 
However$,$ it is also possible to broaden the collection of internal 
semi-norms in various ways. This is done in section 10.4 of reference [8].\par
\vskip 14pt
\noindent {\bf 7. Per-generalized Functions and S-continuity.}
\vskip 14pt
{\bf Theorem 7.1.} {\sl Let generalized function $f \in T$ and standard $p \in \real.$ 
Suppose that $f$ is S-continuous at $p.$ 
Then $f(p) \in {\cal O}.$}\par
\smallskip Proof. Let $f$ be S-continuous at $p \in \real.$ Assume that 
$f(p)$ is not limited. Without loss of generality$,$ assume that $f(p)$ is a 
positive infinite hyperreal number. Now for each $x \approx p,$ it follows 
that $f(x) > (1/2)f(p)$ since by S-continuity $f(x) \approx f(p).$ Note: For
an infinite $\Lambda$  and any $\eps \in \mu(0),$ the infinite 
$\Lambda +\eps > (1/2)\Lambda$.
  By the internal definition method$,$ define that 
internal set 
$$D =\{x\mid(x>0)\land(x \in \hyperreal)\land( \vert x - p \vert 
> 0\to (f(x-p) > (1/2)f(p)))\}.$$
Since $D$ contains all of the positive infinitesimals$,$ then by a modified 10.1.1 in [3], $D$ contains a standard positive $a.$ Consider 
the standard interval $[p-a,p+a].$ It is not difficult to construct a  
 non-negative $h \in {\cal D}$ such that $h(x) = 1,$ for each $x$ such that $\vert x \vert \leq 
(1/2)a$ and $h(x) = 0$ for $\vert x \vert \geq a.$ Now let $g= h(x - p).$
Since $g \in {\cal D},$
$$\hyper {\int_{-\infty}^\infty} f\hyper g = \hyper {\int_{p - a}^{p+a}} 
f\hyper g > 
(1/2)af(p).$$
The results follows from this contradiction.\par
{\bf Corollary 7.1.1.} {\sl Let internal $f \in T$ and standard $p \in \real.$ 
Suppose that $f$ is S-continuous at $p.$  Then any standard function
$F\colon \real \to \real$ such that $F(p) = \st {f(p)}$ is continuous at
$p$ and $f(\monad p) \subset \monad {F(p)}.$} \par
\vskip 14pt
\hrule
\smallskip
\hrule
\smallskip
{\bf Theorem 7.2.} {\sl If $f \in T$ is S-continuous at each $p \in 
\real,$ then 
$F(p) = \st {f(p)}$ is  continuous on $\real.$} \par
\smallskip
\hrule
\smallskip
\hrule
\vskip 14pt
$\{${\bf Remark.} If the function $f\in T$ satisfies the requirements in 
Theorem 7.2$,$ then $f({\cal O})\subset {\cal 
O}.$ If $f$ is also a surjection$,$ then  $f({\cal O})= {\cal 
O}.\}$\par
\vskip 14pt
{\bf Example 7.1} We show that the converse of Theorem 7.2 does not hold.
First$,$ consider for any given standard $a > 0,$ the set $A = [-
a,0)\cup(0,a],$ and the standard function
defined by 
$$f(x) =\cases{1,&$ x \in A$\cr
               0,&$x \in \real - A.$\cr}$$
Then $f^2$ is Riemann integrable on any real $[c,d],\ c \leq d$ and since
$\sup \{f^2(x)\mid x \in \real\} =1,\ \int_c^d f^2 \leq (d -c).$ It follows from 
this that $\vert \int_{-\infty}^\infty fg\vert\leq \vert \int_{-h}^h g\vert\in 
\real$ for some $h \in \real.$ Now consider the internal function $k$
obtained by means of the definition for $f$ but let positive $a = \eps \in 
\monad 0.$ Selecting $c,\ d \in {\cal O},$ it follows that $k \in T.$ Notice 
that for each $p \in \real,\ k(p) = 0.$ Thus $K(p) = \st {k(p)} = 0$ is a 
continuous function on $\real.$ Let $x = \eps/2 \in \monad 0.$ Then 
$k(x) = 1 \not\approx 0.$ Thus $k$ is not S-continuous at $x = 0.$ \par
\vskip 14pt
One of the major advantages in using the nonstandard equivalence class method 
is that various members of a pre-generalized function $\alpha$ can be 
specifically analyzed. This is not the case for  the standard approaches to 
this subject.\par
\vskip 14pt
{\bf Theorem 7.3.} {\sl  If $f \in \alpha$  is S-continuous at each 
$p \in\real,$ then$,$ for the continuous function $F$ defined by 
$F(p) = \st {f(p)}$$,$ $\Hyper F \in\alpha.$ }\par
\smallskip Proof. Since $\Hyper F\in T,$ it follows that $\Hyper F \in \alpha$ if and only if
$\hyper {\int_{-\infty}^\infty} (f - \Hyper F) \hyper g \approx 0$ for all $g \in {\cal D}.$ 
From the definition of $g,$  there is $c \in \real$ such that 
$$\hyper {\int_{-\infty}^\infty} (f - \Hyper F) \hyper g = 
\hyper {\int_{-c}^c} (f - \Hyper F) \hyper g.$$                
Since $[-c,c]$ is compact$,$ if $x \in \Hyper [-c,c],$ then there exists
some $p\in [-c,c]$ such that $x \approx p.$ From S-continuity$,$ standard 
part operator properties$,$ and continuity$,$ $f(x) \approx f(p) \approx 
F(p) \approx \Hyper F(x)$ implies that $f(x) \approx \Hyper F(x).$ 
Consequently$,$ $\vert f(x) - \Hyper F(x) \vert \approx 0.$ Consider the 
internal set $A = \{y \mid (y \in \hyperreal)\land \exists x((x \in \Hyper [-
c,c]) \land(y = \vert f(x) - \Hyper F(x)\vert))\}.$ Since $[-c,c]$ is 
compact$,$ $F([-c,c]) = [a,b]$ implies that $\Hyper F(\Hyper [-c,c]) = \Hyper 
[a,b].$ Further$,$ from above$,$ $f(\Hyper [-c,c]) \subset [-c-1, c + 1].$
Consequently$,$ the internal set  $A$ is *-bounded. Recalling that the
$\Hyper \sup = \sup,$ it follows by *-transfer$,$ that $\sup A \in \hyperreal.$
But for each positive $r \in \real$ and each $y \in A,$ $0\leq y < r.$ Hence,
the $\sup A$ is not a positive real number. Again by *-transfer and the fact 
that $A$ is internal$,$ $\sup A \in \hyperreal.$ Thus $\sup A \approx 0.$    
But$,$ 
$$\left| \hyper {\int_{-c}^c} (f - \Hyper F)   
\hyper g\right| \leq \sup \{\vert f(x) -  
\Hyper F(x)\vert \mid  x \in \Hyper [-c,c]\}\left({\int_{-c}^c}\vert g\vert\right).$$ 
Hence$,$ $\hyper {\int_{-c}^c} (f - \Hyper F)\hyper g \approx 0$ and 
the proof is complete.\par  
{\bf Corollary  7.3.1.} {\sl If $f \in \alpha$ is S-continuous at each $p \in 
\real,$ then $\alpha$ is a ${\cal D}^\prime$-pre-generalized function.}\par
\smallskip Proof. The *-continuous function $\Hyper F \in \alpha.$ Consider a sequence 
$\{g_n\} \subset {\cal D}$ such that $g_n(x) = 0$ for all $n$ and for all $x$ such 
that $\vert x \vert > c \in \real.$ Further suppose that ${g_n} \to 0$ 
uniformly on for all $x$ such that $\vert x \vert \leq c.$ Then 
$\lim_{n \to \infty}\Hyper F[g_n] = \st {\Hyper {(}\lim_{n\to \infty}\langle F, 
 g_n \rangle)} =\st {\hyper {\langle} F, 
\lim_{n\to \infty} g_n \rangle)}=\st 0 = 0.$  \par
\vskip 14pt
The existence of a function $f$ in a pre-generalized function that is 
S-continuous at various members of $\hyperreal$ seems to be of some 
 significance.  Indeed$,$ the standard part of members of a pre-generalized 
function that are 
S-continuous at the same point in $\real$ cannot be distinguished one from 
another at that point.\par
\vskip 14pt
{\bf Theorem 7.4.} {\sl Suppose that $f,h \in \alpha$ and that 
$f$ and $h$ are S-continuous at $p \in \real.$ Then $\st {f(p)}=\st 
{h(p)}$ and if $x \approx p,$ then $f(x) \approx h(x).$ } \par
\smallskip Proof. We know from Theorem 6.1 that $f(p),\ h(p) \in {\cal O}.$ Thus
$\st {f(p)},\ \st {h(p)} \in \real.$ All we need to do is to show that $f(p) 
\approx h(p).$  Since $f,\ h \in \alpha,$ the function $k = f - h \in T_0.$ 
Suppose that $f(p) - h(p)$ is not infinitesimal. Without loss of generality$,$ 
in this case$,$ consider $k(p) =f(p) -h(p) >r >0,\ r \in \hyperreal -\monad 0.$  
As in the proof of Theorem 7.1$,$ there is a $g \in {\cal D},$  such that
$k[g]$ is not infinitesimal; a contradiction.  Thus $f(p) \approx h(p).$ 
Obviously$,$ if $x \in p,$ S-continuity implies that
if $x \approx p,$ then $f(x) \approx f(p) \approx h(p) \approx h(x).$ Hence
$f(x) \approx h(x)$  and the proof is complete.   \par
\vskip 14pt
Due to Theorem 7.4$,$ a pre-generalized function $\alpha$ that contains an 
internal function $f$ that is S-continuous at $p \in \real$ can be considered 
a function itself.\par
\vskip 14pt\vfil\eject
\hrule
\smallskip
{\bf Definition 7.1. ($\alpha$ as a Function on} $\real.$ {\bf )} 
If $f\in \alpha$ and $f$ is S-continuous at $p \in \real$ let
$\alpha (p) = \st {f(p)} = F(p).$ \par
\smallskip
\hrule
\vskip 14pt
Theorem 7.4. shows that definition 7.1 leads to a function-like 
pre-generalized function $\alpha(p).$
However$,$ one other aspect of this $\alpha$ function concept needs to be 
addressed.  How unique is such a $\alpha$ function with respect to members of 
$T$?
\vskip 14pt
{\bf Theorem 7.5.} {\sl Suppose that $f \in \alpha$ is S-continuous
at each $p \in \real.$ 
Let $h \in T$ be S-continuous at each $p \in \real$ and 
$h(p) = f(p).$ Then $h \in \alpha.$}\par
\smallskip Proof. (Notice that simply because $h(p) = f(p)$ at the standard points does 
not imply the functions are equal at the nonstandard points.) What is needed 
is to show that $\langle (f-h), \hyper g\rangle \approx 0$ 
for each $g \in {\cal D}.$ Consider $f-h \in T.$
Let $g \in {\cal D}.$ Then there exists some positive $c \in \real$ such that 
$g(x)= 0 $ for $\vert x \vert \geq c.$ Since $f(\monad a) \subset \monad 
{f(p)}$ and $ h(\monad a) \subset \monad {h(a)}$  and $f(a) = h(a), \ 
(f - h)(\monad a) \subset \monad 0.$ Consider a specific {\it standard} 
positive $\eps$ and the internal set $D(\eps)=$
$$\{y\mid (y \in \hyperreal)\land 
\forall x((x \in \hyperreal)\land(\vert x - a \vert <\eps) \to (\vert f(x) - h(x) \vert < y))\}.$$
Then $D(\eps)$ contains all of the positive infinitesimals. As in the proof of 
Theorem 7.1$,$ there exists a standard positive $\eta_a$ such that such that 
$\vert f(x) - h(x) \vert <\eps$ for each $x$ such that $\vert x - a\vert < 
\eta_a.$ For a fixed $\eps$$,$ for each $a \in [-c,c],$  the positive 
standard  $\eta_a$ exists. The set of open intervals $\{U_a\mid (a \in [-
c,c])\land (U_a = \{ x\mid (a-\eta_a < x < a + \eta_a\})\}$ is an open cover 
of compact $[-c,c]$ and as such there exists a finite subset 
$\{U_1,\ldots,U_n\}$ that covers $[-c,c].$ Also$,$  
$\vert f(x) - h(x) \vert < \eps$ for each $x \in \Hyper U_i.$ Since the set 
$\{U_i\}$ is nonempty and finite it can be adjusted by moving in the end 
points if necessary (remove some of the open intervals) so that 
the shortened open intervals thus obtained have the 
property that no point in $[-c,c]$ is contained in more than two such 
intervals and such that the sum of the lengths is not greater than say $4c.$ 
Denote these adjusted intervals by $H_i.$ Again $\vert f(x) -h(x) \vert < 
\eps$ for each $x \in \Hyper H_i.$ 
This same procedure can be done for any    arbitrary  positive $\eps.$
\par
 Now let 
$\lambda = \max \vert g(x) \vert.$  It is a known fact that $g(x) = g_1(x) +
\cdots + g_n(x),$ where each $g_i \in {\cal D}$ and$,$ for $i = 1,\ldots,n,$ 
$g_i(x) = 0, x \in \real-H_i$ and $\vert g_i(x) \vert \leq \lambda,\ x 
\in \real.$ Hence
$$\left|\hyper {\int_{-\infty}^\infty}(f - h) \hyper g\right| =\left|\sum_{i = 1}^n 
\hyper {\int_{H_i}}(f - h) \hyper g_i \right| \leq$$
$$ \sum_{i = 0}^n \hyper {\int_{H_i}}\eps\lambda \leq 4c\eps\lambda.$$
But$,$ as pointed out$,$ positive $\eps$ is arbitrary. Thus  
$$\left|\hyper {\int_{-\infty}^\infty}(f - h) \hyper g\right|\approx 0$$
and the proof is complete.\par
{\bf Corollary 7.5.1.} {\sl Suppose that $f, \ h\in T$ are S-continuous at each 
$p\in \real$ and $f(p) = h(p)$ for $p\in \real.$ Then there is a unique 
$\alpha$ such that $f, \ h, \ \Hyper F\in \alpha,$ $F = H,$ $F$ is
continuous on $\real$ and  
$\alpha (p) = \st {f(p)} = \st {h(p)}= F(p)$ for each $p \in \real.$}  
\par
\vskip 14pt
Another  aspect of this idea of S-continuity and pre-generalized functions is 
the observation that defining $\alpha (p) = \st {f(p)},$ for some $f \in 
\alpha$ that is S-continuous at $x = p,$ is independent of the 
S-continuous function contained in $\alpha$ by Theorem 7.4. Thus many $\alpha$ 
can be considered as functions on large domains of real numbers. 
For example$,$ the function $d,$ in example 4.1$,$ is in $T$ and is S-continuous 
at every nonzero real number. Thus using this function the Dirac 
${\cal D}^\prime$-pre-generalized function $\delta$ is a function $\delta (x)$ for all nonzero real 
numbers. Then we have  certain algebraic properties associated with the 
$\alpha \in T/T_0$ that are functions at certain standard points. 
\par
\vskip 19pt
\hrule
\smallskip
{\bf Definition 7.2 ($\alpha$ as a Function)} Call a per-generalized function
$\alpha$ a function at  $p \in \real,$ with value $\alpha (p),$ 
if there is an $f \in \alpha$ that is S-continuous at $p$ and let $\alpha (p) 
= \st {f(p)}.$\par 
\smallskip
\hrule
\vskip 14pt
{\bf Theorem 7.6.} {\sl If $\alpha,\ \beta$ are functions at $p,$ then 
$\gamma =\alpha \pm \beta$ is a function at $p$ and $\gamma(p) = 
\alpha(p) \pm \beta(p).$}\par
\smallskip Proof. From the hypothesis$,$ there exist two internal $f,\ h \in T$ such that 
$f, \  h$ are S-continuous at $p$ and $f \in \alpha,\ h\in \beta.$ Hence,
the internal function $f \pm h \in T$ is S-continuous at $p$ and $f \pm h
\in\gamma.$  Since $\st {f(p) \pm h(p)} = \st {f(p)} + \st {h(p)},$ it follows 
that $\gamma (p) = \alpha(p) \pm \beta(p)$ and this completes the proof.
\vskip 14pt
From Theorem 4.3$,$ we know that every member of $\cal F$ is generated by a
member of $\Hyper C^\infty.$ This leads to the concept of the derivative 
of a
generalized function. Unfortunately$,$ if you want to define multiplication for 
the functions in $T$ and use the usual concept that multiplication is 
independent of the member chosen from one or both of two pre-generalized 
functions $\alpha$ 
and $\beta,$  then multiplication must be restricted to certain
per-generalized functions.  Theorem 7.4 states that if $\hyper f,\ \hyper h 
\in \sig C^\infty$ and $\hyper f, \ \hyper h \in \alpha,$ then$,$ since $\hyper f,\ \hyper h $ are S-continuous at 
each $p \in  \real,$ $\st{\hyper f(p)} = f(p) = \st {\hyper h(p)} = h(p).$
Thus $f = h.$ Whenever possible such a unique member of $\sig C^\infty$ can 
be used to generate a unique product-like per-generalized function.
\par
\vskip18pt
{\bf Theorem 7.7.} {\sl Let $\hyper f \in \sig C^\infty$ and $\hyper f\in 
\alpha$ (hence$,$ a ${\cal D}^\prime$-pre-generalized function). 
Suppose that $h,\ k \in \beta.$ Then there exists a pre-generalized 
function $\gamma$ such that $(\hyper f)h,\ (\hyper f)k \in \gamma.$ }\par
\smallskip Proof. Notice that if $\hyper f \in \sig C^\infty$ and $\hyper g \in \sig 
{\cal D},$  
then $\hyper f\hyper g 
\in \sig {\cal D}.$ From Theorem 3.1 (c)$,$ (e)$,$  $(\hyper f)h, (\hyper f)k\in T.$ 
Thus there is a pre-generalized function $\gamma$ such that $(\hyper f)h\in 
\gamma.$ Consider $\langle (\hyper f)h- (\hyper f)k,\hyper g\rangle$ for any 
$g \in {\cal D}.$ Then  $\langle (\hyper f)h- (\hyper f)k,\hyper g\rangle=
\langle h- k,(\hyper f) \hyper g\rangle \approx 0.$ 
\vskip 14pt
\hrule
\smallskip
{\bf Definition 7.3. (Customary Products)} Let $\hyper f\in \sig C^\infty,\ 
\hyper f\in \alpha.$ Then the pre-generalized function
$\gamma$ such that for each $h \in \beta,\ (\hyper f)h \in \gamma$  is the 
customary product of $\alpha$ and $\beta$ and is denoted by $\gamma = 
\alpha\beta.$ \par
\smallskip
\hrule
\vskip 14pt
{\bf Theorem 7.8.} {\sl Let $\{\hyper f,\Hyper h\} \subset \sig C^\infty$ and
$\hyper f \in \alpha,\ \Hyper h \in \beta.$ Then for $\gamma = \alpha\beta$ it 
follows that $\gamma^\prime = \alpha^\prime\beta + \alpha\beta^\prime.$}\par
\smallskip Proof. We know that $(\hyper f\Hyper h)^\prime \in \gamma^\prime\ \hyper f^\prime 
\in \alpha^\prime,\ \Hyper h^\prime \in \beta^\prime.$  Of course,
$(\hyper f\Hyper h)^\prime= (\hyper f^\prime)\Hyper h + (\hyper f)\Hyper 
h^\prime.$ Notice that $(\hyper f^\prime)h \in \alpha^\prime\beta,\ 
(\hyper f)\Hyper h^\prime \in \alpha\beta^\prime.$ Hence$,$  
$(\hyper f^\prime)\Hyper h + (\hyper f)\Hyper 
h^\prime \in \alpha^\prime\beta + \alpha\beta^\prime.$ Then we have that 
$((\hyper f)\hyper h)^\prime \in \gamma^\prime.$ But the pre-generalized 
functions are equivalence classes. Thus $\gamma^\prime = \alpha^\prime\beta + 
\alpha\beta^\prime$ and this completes the proof. \par

{\bf Lemma 7.1.} {\sl Let $f\colon\hyperreal \to \hyperreal$ and 
$h\colon\hyperreal \to \hyperreal$ be S-continuous at $p\in \real.$
Then the product function $k(x) = f(x)h(x)$ is S-continuous at 
$p.$}\par
\smallskip Proof. Let $q \approx p.$ Then  $k(q) -k(p) = f(q)h(q) - f(p)(h(p) = 
f(q) (h(q)-h(p)) + (f(q) - f(p))h(p).$ Notice that $f(q), h(p) \in {\cal O}$ 
and $(h(q)-h(p)), (f(q) - f(p)) \in \monad 0.$ Hence$,$ $k(q) -k(p) \in \monad 
0.$ The result is complete.\par
\vskip 14pt
{\bf Theorem 7.9.} {\sl Let $\hyper f \in \sig C^\infty$ and $\hyper f\in 
\alpha.$ Assume that $\beta$ is a function at $p \in \real.$ Let $\gamma = 
\alpha\beta.$ Then $\gamma$ is a function at $p\in \real$ and $\gamma(p) =
\alpha(p)\beta(p).$}\par
\smallskip Proof. Since $f$ is S-continuous at $p,$ $f(p) = \alpha(p).$ Further$,$ there is 
some $h \in \beta$ that is S-continuous at $p$ and $\beta(p) = \st{h(p)}.$
Consider $k=fh \in \gamma.$ From Lemma 7.1$,$ $\gamma$ is a function at $p$ and
since $\st {f(p)h(p)} = \st {f(p)}\st {h(p)}$ the result follows.\par
\vskip 14pt
Another result similar to Theorem 7.5$,$ but not as definitive$,$ can be obtained 
using the Theorem 7.5 method.\par
\vskip 14pt\vfil\eject
{\bf Theorem 7.10.} {\sl Let $f\in T$ be S-continuous at each $x$
such that $a\leq x\leq b, \ a< b,\ a,b \in \real$ and 
$\alpha$ a function at each $x$
such that $a\leq x\leq b, \ a< b,\ a,b \in \real.$ Then for every
$g \in {\cal D}$ such that the support of $g$ is a subset of $[a,b],$ and every $h 
\in\alpha$ it follows that $f[g]= h[g].$}\par
\smallskip Proof. Let $f \in \beta.$ What is needed is to show that $\langle (f - h), 
\hyper g \rangle \approx 0$ for each $g \in {\cal D}$ such that the support of $g$ is 
a subset of $[a,b].$  Taking the $c$ in the proof of Theorem 7.5 such that
$[-c,c] \subset [a,b],$ the proof is exactly the same as that of Theorem 7.5. 
\par
\vskip 14pt
Any pre-generalized function that is a function on $\real$ is an 
${\cal D}^\prime$-per-generalized function. This idea can be extended to 
the $k$ derivatives of per-generalized functions. \par
\vskip 14pt
{\bf Theorem 7.11.} {\sl Suppose that $\alpha$ is a function on $\real.$
Then for each positive $k \in \nat,$ the pre-generalized function 
$\alpha^{(k)}$ is a ${\cal D}^\prime$-pre-generalized function.}\par
\smallskip Proof. From the hypothesis$,$ there is a function $F$ continuous on $\real$ such 
that $\Hyper F\in \alpha.$ From Theorem 4.3$,$ there is some $h \in \Hyper
{C^\infty}$ such that $h\in \alpha$ and$,$ from the Definition of $\alpha^{(k)},$ 
$h^k \in \alpha^{(k)}.$ Let $\{g_n\}$ be a sequence of members of ${\cal D}$ such 
that $g_n(x) = 0$ for all $n$ and all $x$ such that $\vert x \vert
>c, \ c > 0$ and $\{g_n^{(k)}\}$ converges uniformly in $x$ for all positive
$k.$ All that is needed is to show that the sequence $h^{(k)}[g_n]$ converges 
to zero. We know by parts integration that $h^{(k)}[g_n]= (-1)^k h[g_n^{(k)}].$
Further$,$ $\Hyper F[g_n^{(k)}]= h[g_n^{(k)}].$ From uniform convergence,
$\lim_{n \to \infty} \Hyper F[g_n^{(k)}] = \lim_{n \to \infty} h[g_n^{(k)}]=
\lim_{n \to \infty} h^{(k)}[g_n]=  0.$ The result follows.\par
\vskip 14pt
\noindent {\bf 8. Generalizations and Beyond} \par
\vskip 14pt
Every theorem in the previous sections relative to pre-generalized functions,
holds true if the domain $\hyperreal$ is changed to an open domain  $\Omega 
\subset \hyperrealp m.$
[Note: Theorem 4.3 still holds in the sense that there exists $f\in \Hyper 
C^\infty(\real^m)$ that generates each pre-generalized function.]  
Depending upon the level of pre-generalized function differentiation desired$,$ 
the test space ${\cal D}$ can be enlarged or reduced. For example$,$ if you only 
wish that 
for each $\alpha,$  there exists $\alpha^{(k)}$ where $1 \leq k \leq m.$
Then the test spaces can be the set of all functions $g \colon \real \to 
\real$ that have bounded support and are continuously differentiable of order
$m.$ In this case$,$ all the previous theorems on pre-generalized functions$,$ 
modified for this degree of differentiability$,$ also hold. Then another 
approach that has proved  to be accessible to nonstandard methods is the 
generalization due to Colombeau [8]. I mention that the new approach by 
Egorov [9] could lead to a very accessible nonstandard theory and needs to be 
fully explored.\par
Using section 6 and an appropriate collection of internal semi-norms,
 the concepts of the S-limit and S-convergence$,$ and the like$,$ 
can now be applied to  
pre-generalized functions.   
\par
\vskip 14pt
\centerline{\bf References}
\bigskip
\id{\bf 1.} Albeverio and et el. 1986. {\it Nonstandard Methods in Stochastic 
Analysis and Mathematical Physics,} Academic Press$,$ Orlando.
\id{\bf 2.} Davis$,$ M. 1977. {\it Applied Nonstandard Analysis,} John Wiley \& Sons$,$ 
New York.
\id {\bf 3.} Herrmann$,$ Robert A. 1991. Nonstandard Analysis Applied to 
Advanced Undergraduate Mathematics - Infinitesimal Modeling and {\it 
Instructional Development Project$,$ Mathematics Department$,$ U. S. Naval 
Academy$,$ Annapolis$,$ MD 21402-5002.}\hfil\break 
http://www.arXiv.org/abs/math.GM/0312432
\id{\bf 4.}{\vbox{\hrule width 0.75in}} 1989. Fractals and ultrasmooth microeffects. {\it 
J. Math. Physics}. 30:805--808. 
\id{\bf 5.} Hurd$,$ A.E. and P. A. Loeb. 1985. {\it Introduction to Nonstandard 
Real Analysis,} Academic Press$,$ Orlando.
\id{\bf 6.} McShane$,$ Edward J. 1947. {\it Integration,} Princeton University 
Press$,$ Princeton.
\id{\bf 7.} Rudin$,$ Walter. 1953. {\it Principles of Mathematical Analysis,} 
McGraw-Hill$,$ New York.
\id{\bf 8.} Stroyan$,$ K. D and W. A. J. Luxemburg. 1973 {\it Introduction to 
the Theory of Infinitesimals,} Academic Press$,$ New York.\par
\bigskip
\centerline {\bf Additional References}\vskip 14pt
\id{\bf 9.} Colombeau$,$ J. F. 1992. {\it Multiplication of Distributions,} 
Lecture Notes in Mathematics$,$ \#1532$,$ Springer-Verlag$,$ New York.
\id{\bf 10.} Egorov$,$ Y. V. 1990. A contribution to the theory of generalized 
functions$,$ {\it Russian Math. Surveys,} 45(5): 1--49.

\end